\documentclass[11pt]{article}
\usepackage{amssymb,amsfonts,latexsym,a4}
\usepackage{graphics,color}
\usepackage{epsfig, subfigure}

\newcommand{\lastcorrections}%
{{
 \begin{sloppypar}
    \baselineskip -0.2in
    \tiny\bf\noindent
MJG:
Version of
Tue Jun 28 00:09:09 HKT 2005

\end{sloppypar}
}}

\newenvironment{Proof}%
        {\par\noindent{\bf Proof.}}%
        {{\hfill{\Large$\Box$}}\smallskip}
\newtheorem{Lemma}{Lemma}

\newtheorem{Cor}{Corollary}
\newtheorem{Definition}{Definition}

\newcommand{\Cn}{{C_n}}
\newcommand{\Ln}{{L_n}}
\newcommand{\Lnpo}{{L_{n+1}}}

\newcommand{\Hook}{{\mbox{\tt Hook}}}
\newcommand{\Hookn}{{\mbox{\tt Hook}(n)}}
\newcommand{\N}{{\mbox{\tt New}}}
\newcommand{\Nn}{{\mbox{\tt New}(n)}}

\newcommand{\ID}[2]{{\mbox{\tt ID}_{#1}(#2)}}
\newcommand{\OD}[2]{{\mbox{\tt OD}_{#1}(#2)}}

\newcommand{\perm}{{\mbox{\rm Perm}}}

\newcommand{\calp}{{\cal P}}

\newcommand{\call}{{\cal L}}

\newcommand{\calcc}{{\cal CC}}

\newcommand{\Llt}{{<_L}}
\newcommand{\Rlt}{{<_R}}
\newcommand{\Leq}{{=_L}}

\newcommand{\Lleq}{{\leq_L}}
\newcommand{\Rleq}{{\leq_R}}
\newcommand{\ems}{{\bar s}}

\title{Permanents of Circulants:\\ a Transfer Matrix Approach\\ 
(Expanded Version)\footnote{\small Partially
supported by HK CERG grants HKUST6162/00E and 613105.
Dept. of Computer Science \& Engineering, Hong Kong U.S.T., Clear Water Bay,
Kowloon, Hong Kong.  Email addresses are {\em \{golin,cscho,yalding\}@cse.ust.hk}.
}
}
\author{Mordecai J. Golin
\and         Yiu Cho Leung
\and Yajun Wang
}

\date{}
\begin{document}
\thispagestyle{empty}
\maketitle
\begin{abstract}
Calculating the permanent of a $(0,1)$ matrix is a $\#P$-complete problem but there are some classes
of {\em structured} matrices for which the permanent is calculable in polynomial time.  The most
well-known example is the {\em fixed-jump} $(0,1)$
circulant matrix  which, using algebraic techniques,
was shown by Minc to satisfy a constant-coefficient fixed-order recurrence relation.

In this note we show how, by interpreting the problem as calculating the number of 
cycle-covers in a directed circulant {\em graph,} it is straightforward to reprove Minc's
result using  combinatorial methods.  This is a two step
process: the first step is to show that the cycle-covers of directed circulant
graphs  can be evaluated using a 
{\em transfer matrix} argument.  The second is to show that the associated transfer matrices,
while very large, actually have much smaller characteristic polynomials than would a-priori be expected.
 
An important consequence of this new viewpoint is that, in combination
with a new
recursive decomposition of circulant-graphs,  it permits extending
Minc's result to calculating the permanent  of the much larger
class of circulant matrices with
{\em non-fixed} (but linear) jumps.  It also permits us to count 
other types of structures in circulant graphs,  e.g., Hamiltonian Cycles.

\end{abstract}

\thispagestyle{empty}

\section{Introduction}
\label{sec:Introduction}
\begin{Definition}
Let $A=(a_{i,j})$ be an $n \times n$ matrix.
Let $S_n$ be the set of permutations of the integers
$[1,\ldots,n]$.  The {\em permanent} of $A$ is
\begin{equation}
\label{eq:permdef}
\perm(A)= \sum_{\pi \in S_n} \prod_{i=1}^n a_{i,\pi(i)}
\quad\mbox{where}\quad
\pi = [\pi(1),\ldots,\pi(n)].
\end{equation}
\end{Definition}

If $A$ is a $(0,1)$ matrix, then  $A$ can be interpreted as the adjacency matrix
of some directed graph $G$ and $\perm(A)$ is the number of {\em directed cycle-covers} in $G$,
where a directed cycle-cover is a collection of disjoint cycles that cover all of the
vertices in the graph.
Alternatively, $A$ can be interpreted as the adjacency matrix of
a bipartite graph $\bar G$, in which case 
$\perm(A)$ is the number of {\em perfect-matchings} in $\bar G.$
The permanent is a classic well-studied combinatorial object (see the book and
later survey by Minc\cite{Minc78,Minc87b}).

Calculating the permanent of a $(0,1)$ matrix is a $\#P$-Complete problem \cite{Valiant79}
even when $A$ is restricted to have only 3 non-zero entries per row \cite{DaLuMiVa88}.
The best known algorithm for calculating a general permanent is a straightforward
inclusion-exclusion  technique due to Ryser \cite{Minc78} running in $\Theta(n 2^n)$ time 
and polynomial space.
By allowing super-polynomial space, Bax and Franklin \cite{BaFr02} developed
a slightly faster (although still exponential) algorithm for the $(0,1)$ case. 
For non-exact calculation Jerrum, Sinclair and  Vigoda  
\cite{JeSiVi04}
have 
developed a fully polynomial approximation scheme for {\em approximating} the
permanent of nonnegative matrices. 

On the other hand, for certain special structured 
classes of matrices one can {\em exactly} calculate the 
permanent in ``polynomial time''. The most studied example of such a class is probably 
the {\em circulant matrices,} which, as discussed in \cite{CoShWi02}, can be thought of as
the borderline between the easy and hard cases.

\begin{figure*}[t]
$$
\begin{array}{c}
  \subfigure[$P_n^{-1}+ I + P_n^2$ ($n=6$)]
  {\parbox[b][1.2in][t]{1.8in}
    {
    \vspace{-.1in}
    \begin{center}$
    \left(
       \begin{array}{cccccc}
          1 & 0 & 1 & 0 & 0 & 1 \\
          1 & 1 & 0 & 1 & 0 & 0 \\
          0 & 1 & 1 & 0 & 1 & 0 \\
          0 & 0 & 1 & 1 & 0 & 1 \\
          1 & 0 & 0 & 1 & 1 & 0 \\
          0 & 1 & 0 & 0 & 1 & 1
       \end{array}
    \right)
    $\end{center}
    }
  }
  \hspace{.6in}
  \subfigure[$C_{n}^{-1,0,2}$ ($n=6$)]{\epsfig{file=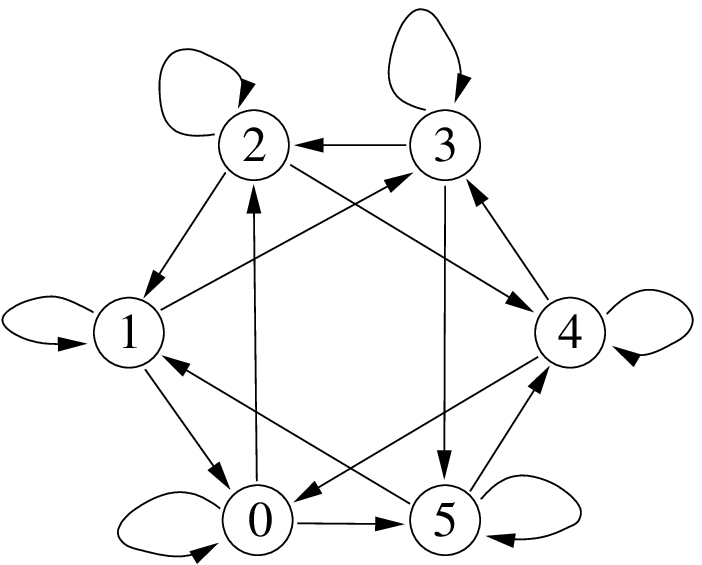, width= 1.4 in}} \\
  \subfigure[$P_n^{-1}+ I + P_n^2$ ($n=7$)]
  {\parbox[b][1.4in][t]{1.8in}
    {
    \vspace{-.1in}
    \begin{center}$
    \left(
       \begin{array}{ccccccc}
          1 & 0 & 1 & 0 & 0 & 0 & 1 \\
          1 & 1 & 0 & 1 & 0 & 0 & 0 \\
          0 & 1 & 1 & 0 & 1 & 0 & 0 \\
          0 & 0 & 1 & 1 & 0 & 1 & 0 \\
          0 & 0 & 0 & 1 & 1 & 0 & 1 \\
          1 & 0 & 0 & 0 & 1 & 1 & 0 \\
          0 & 1 & 0 & 0 & 0 & 1 & 1
       \end{array}
    \right)
    $\end{center}
    }
  }
  \hspace{.4in}
  \subfigure[$C_{n}^{-1,0,2}$ ($n=7$)]{\epsfig{file=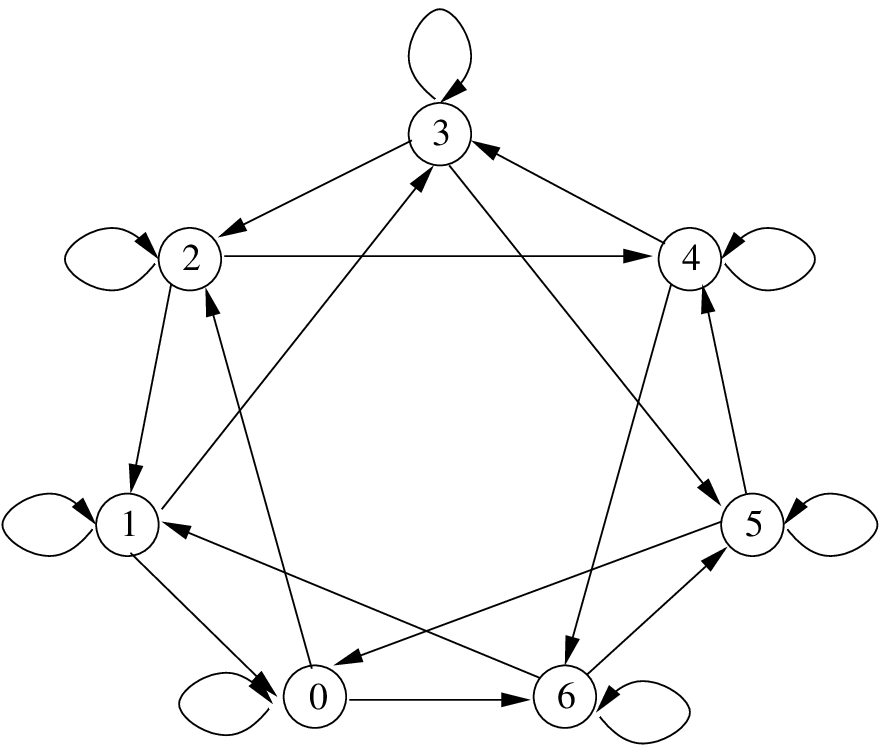, width= 1.8 in}} 
\end{array}
$$

\caption{$C_{n}^{-1,0,2}$: Circulant matrices (a) and (c) are, respectively,
 the adjacency matrices of circulant graphs
$C_{n}^{-1,0,2}$ in (b) and (d) for $n=6, 7$.
}
\label{fig:const}
\end{figure*}

An $n \times n$ circulant matrix $A=(a_{i,j})$ (see Figures \ref{fig:const} (a) and (c))
is defined by specifying
its first row;  the
$(i+1)^{\mbox{st}}$ row is a cyclic shift $i$ units to the right of the first row, i.e.,
$a_{i,j}= a_{1,1+(n+j-i)\bmod n}.$
Let $P_n$ denote the $(0,1)$ $n \times n$ matrix with $\bf 1$s in positions
$(i,i+1)$, $i=1,\ldots,n-1,$ and $(n,1)$ and $\bf 0$s everywhere else.
Many of the early  papers on this topic express circulant matrices in the form
\label{eq:A_def}
\begin{equation}
A_n=a_1 P_n^{s_1} + a_2 P_n^{s_2} + \cdots + a_k P_n^{s_k}
\end{equation}
where 
$0 \le s_1 < s_2 < \cdots < s_k < n$
and $ a_i = a_{1,{s_i+1}}.$

The first major result on permanents of $(0,1)$ circulants was due to Metropolis, Stein and Stein \cite{MeStSt69}. 
Let $k>0$ be fixed and $A_{n,k} = \sum_{i=0}^{k-1} P^i_n,$ be the $n \times n$ circulant matrix whose first
row is composed of $\bf 1$s in  its first $k$ columns and $\bf 0$s  everywhere else.  Then
\cite{MeStSt69} showed that, as a function of $n$, $\perm(A_{n,k})$ satisfies a fixed order constant-coefficient
recurrence relation in $n$ and therefore, could be calculated in polynomial time in $n$ (after a superpolynomial
``start-up cost'' in $k$ for deriving the recurrence relation).

This result  was greatly improved by Minc who showed that it was only 
a very special case of a general rule. 
Let  $0 \le s_1 < s_2 < \cdots < s_k < n$ be {\em any} fixed sequence and set 
$A_n = A_n(s_1,\ldots,s_k) =  P_n^{s_1} +  P_n^{s_2} + \cdots +  P_n^{s_k}$. In \cite{Minc85,Minc87a} Minc
proved that  
{\em $\perm(A_n)$}
always satisfies a constant-coefficient
recurrence relation in $n$ of order $2^{s_k}-1.$ Minc's theorem was proven by
manipulating algebraic properties of $A_n.$ 
Note, that as mentioned by Minc, this result is difficult to apply for large $s_k$ since, in order
to derive the coefficients of the recurrence relation it is first necessary to evaluate
$\perm(A_n)$
for $n\le 2(2^{s_k}-1)$ and, using Ryser's algorithm,  this requires
$\Omega\left(2^{2^{s_k}}\right)$ time.

Later   Codenotti, Resta  and various coauthors improved these results in various ways; e.g. 
in \cite{BeCoCrRe99} showing how to evaluate {\em sparse} circulant matrices of size $\le 200$;
in \cite{CoCrRe97, CoRe01} showing that the permanents of circulants with only three $\bf 1$s per
row can be evaluated in polynomial time; 
in \cite{CoRe02} showing how the permanents of some special sparse circulants can be 
expressed in terms of determinants and are therefore solvable in polynomial time; 
in \cite{BeCoCrRe99} showing that the permanents of {\em dense} circulants are hard to calculate
and in \cite{CoShWi02} that even approximating the permanent of an arbitrary circulant 
modulo a prime $p$ is ``hard'' unless $\mbox{\bf  P}^{\mbox{\bf \#P}} = \mbox{\bf BPP}$.

In this paper we return to the original problem of Minc.
Our first main result  will be to show that if circulant {\em matrix}
$A_n(s_1,\ldots,s_k)$ is interpreted as the adjacency matrix of  a directed circulant
{\em graph}  $\Cn$, then   counting the number of
cycle-covers of $C_n$ using  a {\em transfer matrix} approach  immediately reproves 
Minc's result. In addition to  rederiving Minc's original result using  a combinatorial rather
than algebraic proof  this new technique
permits  us extend the result to a much larger set of circulant graphs.  It will
also permit  us to  
address other problems, e.g., counting Hamiltonian cycles in circulant graphs,
which at first might seem unrelated. To explain, we first need to introduce some notation.
\begin{Definition} See Figure \ref{fig:const}.
\label{def:circ_const}
Let  $C_n^{s_1,s_2,\cdots,s_k}.$ be 
the  $n$-node {\em directed  circulant graph} with jumps $S = \{s_1,s_2,\ldots s_k\}$.
(Note that this definition permits  negative $s_i$.)
Formally, 
$$C^{s_1,s_2,\dots,s_k}_n= (V(n),E_C(n))$$ 
where 
$$V(n)=\{0,1,\dots,n-1\}$$
and 
$$E_C(n)=\Big\{(i,j) : (j-i) \bmod n \in S\Big\}.
$$
{\small \em Note: we will assume that $S$ contains at least
one non-negative $s_i$ since, if all the $s_i$ were negative, 
we could multiply them by $-1$ and get
an isomorphic graph.  Also,  we will often write $\Cn$ as shorthand for  $C_n^{s_1,s_2,\cdots,s_k}.$
}
\end{Definition}

Let $G=(V,E)$ be a graph, $T \subseteq E$ and $v \in V.$ Define $\ID T v$ to be the {\em indegree} of 
$v$ in graph $(V,T)$ and $\OD T v$ to be the {\em outdegree} of 
$v$ in  $(V,T)$.
$T  \subseteq E$ is a {\em cycle-cover} of $G$ if
\begin{equation}
\label{eq:cycle_condition}
\forall v \in V,\quad \ID T v = \OD T v = 1.
\end{equation}

\begin{Definition}
Let $S = \{s_1,s_2,\ldots s_k\}$ be given.
Set 
$${\cal CC}(n) = \{T \subseteq C_n \,:\, \mbox{ $T$ is a cycle-cover of $C_n$}\}$$
and
$$T(n) = |{\cal CC}(n)| = \mbox{No. of cycle-covers of $C_n$}.$$
\end{Definition}

Note that, by the standard correspondence mentioned previously,
$A_n(s_1,\ldots,s_k)$ is the adjacency matrix of $C_n^{s_1,s_2,\cdots,s_k}$
and $T(n) = \perm(A_n(s_1,\ldots,s_k))$.  So,  calculating $T(n)$ is 
equivalent to calculating permanents
of $A_n(s_1,\ldots,s_k)$.

There is also a well-known simple correspondence between cycle covers and permutations.
Consider the {\em directed} complete graph with all $n^2$ distinct edges on $n$ vertices (self-loops are
permitted). Now let 
${S}_n$ be the set of $n!$ permutations on $[0,\ldots,n-1]$.  For
a fixed permutation $\pi \in S_n$, the set of edges $\bigcup_{i=0}^{n-1} (i,\pi(i))$
is a cycle cover.  In the other direction suppose $T$ is  a cycle cover.  Define $\pi$
by 
 $\pi(i) =j$ where $j$ is the unique vertex such that $(i,j) \in T$.  Then $\pi$
 is a permutation.
This is a one-one correspondence between cycle covers and permutations so $T(n)$ 
counts the number of permutations $\pi \in S_n$  restricted such  that 
 $(\pi(i)-i) \bmod n \in S.$  For example,  if $S=\{1,2,3\}$, the number of cycle covers in
 the corresponding circulant graph $C_n^{1,2,3}$ is equal to
the number of permutations $\pi$
 such that  $\Bigl( (\pi(i)-i) \bmod n \Bigr) \in \{1,2,3\}.$
 In fact, in \cite[Sec 4.7]{Stanley86}, Stanley 
shows that, for fixed $S$,  the number of such permutations satisfies a recurrence relation,
giving an alternative derivation of Minc's result for this special case
(but without the bound on the
order of the recurrence relation given in \cite{Minc85,Minc87a}).

In \cite{GoLe04,GoLeWa04} the authors of this paper 
were interested in counting spanning trees and other structures in
{\em undirected circulant} graphs.  The main tool introduced there was
a recursive decomposition of such graphs.  In Section \ref{sec:decomposition}
we describe a related recursive decomposition of {\em directed circulant} graphs.
Our technique will  be to use this decomposition to show that 
for some constant $m$ there is a 
$m \times 1$ (column) vector function $\bar T (n)$ such that 
\begin{equation}
\label{eq:trans_intro}
\hspace*{.15in} 
\forall n \ge 2\ems,\quad 
T(n) = \beta\,   \bar T (n)
\quad \mbox{and} \quad 
\bar T(n+1) =A \, \bar T(n)
\end{equation}
where  $\ems$ is a constant to be defined later (but  reduces to $\bar s = s_k$ 
for the Minc formulation
described previously), $\beta$ is a $1 \times m$ constant row-vector and 
$A$ is a constant $m \times m$ matrix. Such an $A$ is known as a {\em transfer-matrix} see, e.g., 
\cite{Stanley86}.
 
Let $P(x)=\sum_{i=0}^t p_i x^i$ be any polynomial that  annihilates $A,$ i.e., 
$P(A) =0$.   Then it is easy to see that $\forall n \ge 2\ems,$
\begin{eqnarray*}
\sum_{i=0}^t p_i T(n+i)
&=& \beta  \left(\sum_{i=0}^t p_i A^{n+i-2 \ems}\right)  \bar T (2\ems)\\
&=& \beta\,  A^{n-2\ems} \left(\sum_{i=0}^t p_i A^{i}\right)  \bar T (2\ems)\\
&=& \beta\,   A^{n-2\ems} \,  {\bf 0}\,   \bar T (2\ems) \\
&=& 0
\end{eqnarray*}
where ${\bf 0}$ denotes the $m \times m$ zero matrix and $0$ a scalar; $T(n)$ thus satisfies the
degree-$t$  constant coefficient recurrence relation 
$T(n+t) = \sum_{i=0}^{t-1}  -\frac {p_i} {p_t} T(n+i)$ 
in $n$. 
By the Cayley-Hamilton theorem, the characteristic polynomial of $A$
-- which has degree
$\le m$ --
must annihilate
$A,$ so such a polynomial exists and $T(n)$ satisfies a recurrence relation of at most
degree $m.$
In our notation, Minc's theorem is that $T(n)$ satisfies a recurrence relation of degree $2^{\ems} -1.$ 
Unfortunately, in our construction, $m = 2^{2 \ems}$ so the characteristic
polynomial does not suffice for our purposes. Our next step 
will involve showing that even though $A$ is of size  $2^{2  \ems}\times 2^{2  \ems}$, 
there is a much smaller  $P$, of degree
$2^{\ems}-1$, that annihilates $A$, thus reproving Minc's theorem. 
We point out that this degree reduction of the transfer matrix (to less
than the square-root of the original size)
is, a-priori,  quite unexpected, and does not occur in the undirected-circulant counting  problems analyzed in
\cite{GoLe04,GoLeWa04}.

One interesting consequence of this new derivation is that, unlike in Minc's proof,  to derive the recurrence
relation it is no longer necessary to start by spending
$\Omega\left(2^{2^{\ems}}\right)$ time calculating the first $2^{\ems}$ 
values of $T(n)$ using Ryser's method.  Instead one only has to calculate 
$A$, $\beta$, the polynomial $P$ and the first $2^{\ems}$ values of $\bar T(n)$ which, as we will see later,
 can all be done in 
$O(\ems 2^{4\ems })$ time, reducing the start-up complexity from  doubly-exponential in 
$\ems$ to singularly exponential. 

Another, albeit minor, consequence of this new derivation
is that it can also handle non-$(0,1)$ circulants. That is,  given {\em any} matrix $A_n$ of
the form (\ref{eq:A_def}), even when the $a_i$ are not restricted to be
in $\{0,1\}$  the technique shows that 
$\perm(A_n)$  satisfies a recurrence relation of degree $2^{\ems} -1.$
This is only a minor consequence, though,
since  working through the details of Minc's original proof it is possible to modify
it to get the same result.


\begin{figure*}[h]
$$
\begin{array}{c}
\subfigure[$C_{3n}^{1,n,2n}$ ($n=4$)]{\epsfig{file=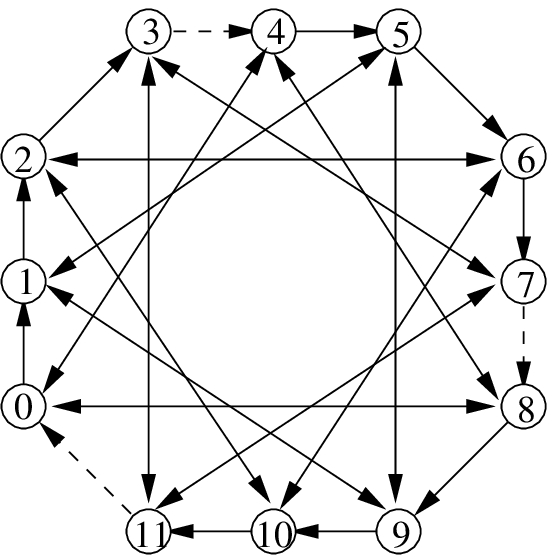, width= 1.4 in}}
  \hspace{.6in}
\subfigure[$C_{3n}^{1,n,2n}$ ($n=5$)]{\epsfig{file=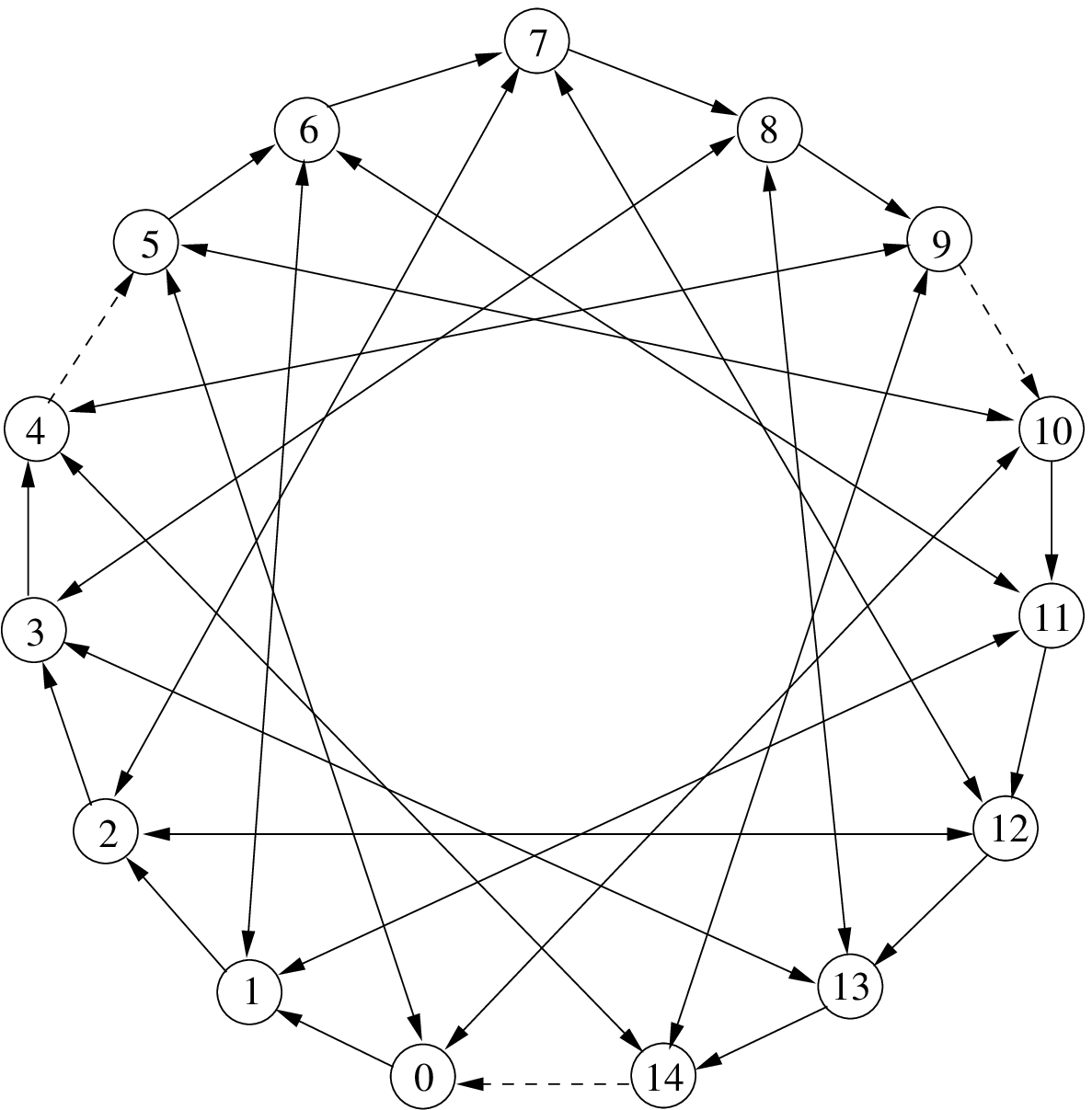, width= 2.0 in}}\\
\subfigure[$C_{3n}^{1,n,2n}$ ($n=4$)]{\epsfig{file=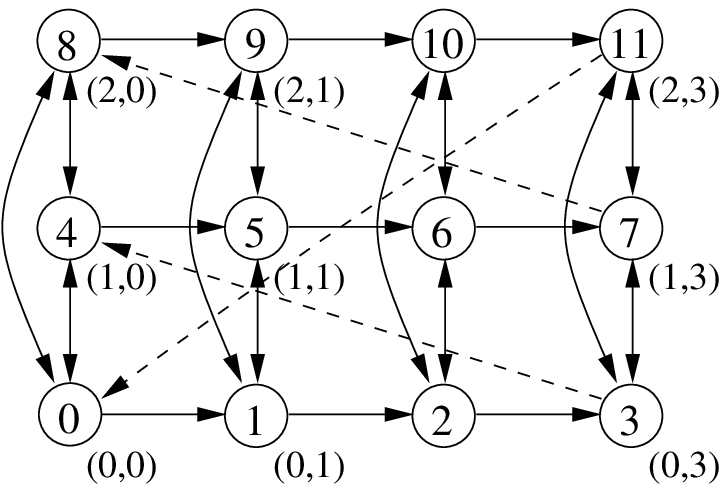, width=2.1 in}}
  \hspace{.2in}
\subfigure[$C_{3n}^{1,n,2n}$ ($n=5$)]{\epsfig{file=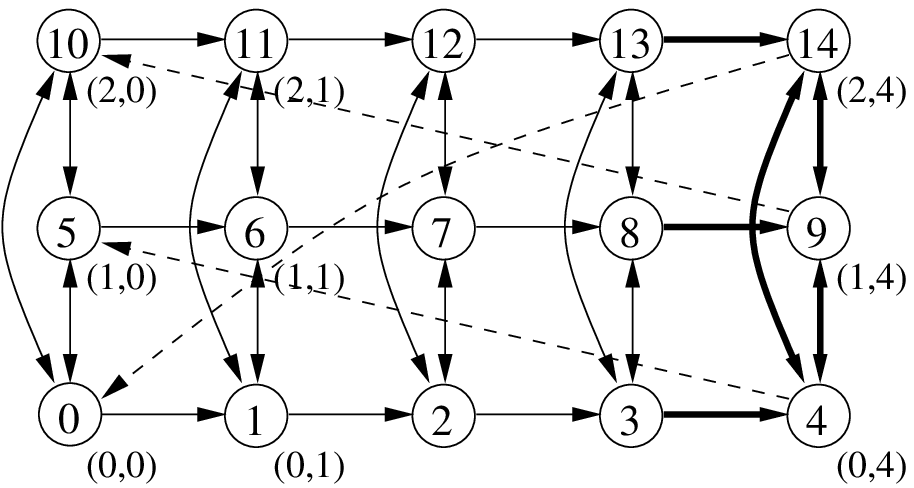, width=2.65 in}}
\end{array}
$$
\caption{$C_{3n}^{1,n,2n}$, a non-constant jump circulant: 
Solid edges are $L_n.$
Dashed edges are $\Hookn.$ 
(a) and (b) are the circulant graphs when $n=4,5$.
(c) and (d) are corresponding lattice representations of the same graphs.
The bold solid  edges on the right of (d) are $\N(5)=L_5-L_4$.
The 3 vertices$\{4,9,14\}$  on the right are $VN(n).$
Note that the dashed $\Hookn$ edges for both $n=4,5$ are ``independent'' of $n.$}
\label{fig:non_const}
\end{figure*}

A much  more important new  consequence, and a major motivation for this paper,
is the fact that  the proof can be extended to evaluate the
permanents of 
{\em non-constant (linear)} 
jump circulant matrices, something which has not been addressed before.
As an example Minc's technique would not permit calculating the permanents of
$A_{3n}(1,n,2n)$, something which our new method allows.
To explain this, we generalize Definition \ref{def:circ_const} to
\begin{Definition} See Figures \ref{fig:non_const} (a) and (b).
\label{def:circ_non_const}
Let $p,s$,    $p_1,p_2,\ldots,p_k$  and $s_1,s_2,\ldots,s_k$
be fixed   integral constants with  such that 
$\forall i,$ $0 \le p_i <p.$ Set 
$S = \{p_1n+s_1,p_2n+s_2,\cdots,p_kn+s_k\}.$
Denote the  $(pn+s)$-node {\em directed  circulant graph} with jumps $S$ 
 by  
$$\Cn = C_{pn+s}^{p_1n+s_1,p_2n+s_2,\cdots,p_kn+s_k} = (V(n),E_C(n))$$
where 
$$V(n)=\{0,1,\dots,pn+s-1\}$$
and
$$E_C(n)=\Big\{(i,j) : (j-i) \bmod (pn+s) \in S\Big\}.
$$
\end{Definition}

Figure (2a) and (2b) illustrate $C^{1,n,2n}_{3n}$ for $n=4,5$. Figure (2c) and (2d) are the
corresponding lattice representation, which will be introduced in section 4.

Note that $A_{pn+s}(p_1n+s_1,p_2n+s_2,\cdots,p_kn+s_k)$ is the adjacency
matrix of $\Cn$ so, counting the cycle-covers in $\Cn$ is equivalent
to evaluating
$\perm(A_{pn+s}(p_1n+s_1,p_2n+s_2,\cdots,p_kn+s_k))$. 
Our method of counting the cycle covers in $\Cn$ will be to derive a new recursive
decomposition of $\Cn$ (which might be of independent interest) and use
it to show that an analogue of (\ref{eq:trans_intro}) holds in the non-constant jump
case as well; thus  $T(n)$ still 
satisfies a constant-coefficient
recurrence relation in $n$.
For example,  Table \ref{tab:cov_results}, shows the recurrence
relation for the number of cycle covers in $C^{1,n+1,2n}_{3n}$, $C^{0,n,2n-1}_{3n}$
, $C^{1,n,2n+1}_{3n+1}$ and $C^{2,n+1, 2n+2}_{3n+1}.$


In the next section
 we describe the  recursive decomposition of $C_n$, for  constant-jump circulants upon
which our technique is based.  In Section \ref {sec:Minc} we show how this permits easily reproving
Minc's result for constant-jump  circulants.  
In Section \ref{sec:non_constant} we then describe the generalization of the decomposition
and the  minor modifications to the proofs
that are needed to extend our analysis to the non-constant circulants introduced in
Definition \ref{def:circ_non_const}.
Finally, in Section \ref{sec:Extensions}, we sketch generalizations and
other uses of our technique;  we first show how it can be extended to calculate permanents
of non 0-1 circulants.  We then describe how it can be used to calculate the moments
of the the random variable counting the number of cycles in a random  restricted permutations.
We conclude by discussing   how to extend the technique to counting the number of Hamiltonian cycles in
directed circulants, extending the result of \cite{YaBuCeWo97},
 which only worked for circulant graphs with two jumps.



\begin{table*}[t]
$$
\begin{array}{|l|l|c|}
\hline
C_n^{-1,0,1}  & T(n) = 2T(n-1) - T(n-3) & T(n) \sim \phi^n \\
C_n^{0,1,2}   & \mbox{ initial values } 9,13,12 \mbox{ for } n=4,5,6 & \phi = (1 + \sqrt 5)/2 \\
\hline
 C^{0,n,2n-1}_{3n} &  \begin{array}{l} T(n) = 5T(n-1) - 5T(n-2) \\
  \hspace{.55in} -5T(n-3) + 6T(n-4) \end{array} & \\  
 C^{1,n+1,2n}_{3n}& \mbox{ initial values } 17,45,113,309 \mbox{ for }
n=2,3,4,5 & T(n) \sim 3^n\\
\hline
 C^{1,n,2n+1}_{3n+1} & T(n)=4T(n-1)+5T(n-2) & \\
                     & \hspace{.55in} -16T(n-3)-2T(n-4) & \\
C^{2,n+1,2n+2}_{3n+1} & \hspace{.55in} -8T(n-5)-6T(n-6) & \\
                      & \hspace{.55in} +16T(n-7)+3T(n-8) & T(n) \sim \psi\phi^n\\
                      & \hspace{.55in} +4T(n-9)+T(n-10) & \psi = (1 + \sqrt 5)/2 \\
& \mbox{ initial values } 31,169,523,2401,9351,40401, & \phi = 2 + \sqrt 5\\
& \mbox{ } 167763,714025,3010351,12766329 & \\
& \mbox{ for } n=2,3,\dots,11 & \\
\hline
\end{array}
$$
\label{tab:cov_results}
\caption{The number of cycle-covers $T(n)$ in directed circulant graphs with constant jumps
$C_n^{-1,0,1}$ and $C_n^{0,1,2}$,  
and with non-constant jumps   $C^{0,n,2n-1}_{3n},$ and $C_{3n}^{1,n+1,2n}$, 
and  $C_{3n+1}^{1,n,2n+1}$ and $C_{3n+1}^{2,n+1,2n+2}$
as derived by the techniques of this paper. 
Note that for all pairs of graphs, the number of cycle covers for each of the graphs in
the pair  is  the same.  
This is because the adjacency matrices for the two items in each pair are just linear circular shifts of each other so the permanents
of their adjacency matrices are the same.  The second item in each pair is in the form that
we analyze.  That is,  for the constant case,  having $s_1=0$, and for the nonconstant case,
having, $\forall i,$  $s_i \ge s.$}
\end{table*}




\section{A Recursive Decomposition of Directed Circulant Graphs}
\label{sec:decomposition}

The main conceptual difficulty with deriving a recurrence relation for 
$T(\Cn)$ is that larger circulant graphs can not be built recursively
out of smaller ones.  The crucial observation,  though, is that, there is
{\em another graph}, $\Ln$, the {\em lattice   graph}, that {\em can} be built recursively, and
$\Cn$ can then be constructed from $\Ln$ through the addition of a constant number of 
edges\footnote{To put this into context, this is very similar to the definition of 
{\em Recursive families} for undirected graphs \cite{BiDaSa72,NoRi04},
which were used for recursively building the Tutte polynomials of graphs in a class.}. In \cite{GoLe04,GoLeWa04}
the authors of this paper developed such a recursive
decomposition for {\em undirected} circulant graphs as a tool for counting the number of spanning trees
in such graphs.   In what follows we develop a corresponding decomposition for
{\em directed circulants} that will permit counting cycle-covers. 

We first show this for the restricted case in which $S$, the set of jumps, is constant (independent of $n$), where
it is easy to visualize.  In Section \ref{sec:non_constant}
we will see how to  extend the decomposition to the more complicated
case in which the set of jumps can depend linearly upon $n$, as described in Definition \ref{def:circ_const}.

\begin{figure*}[t]
$$
\begin{array}{c}
\subfigure[$C_{6}^{0,1,2}$]{\epsfig{file=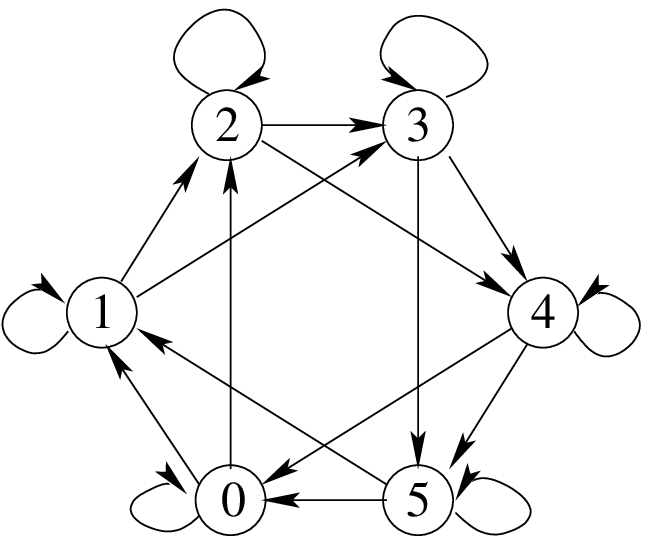, width= 2.0 in}}
  \hspace{.4in}
\subfigure[$C_{7}^{0,1,2}$]{\epsfig{file=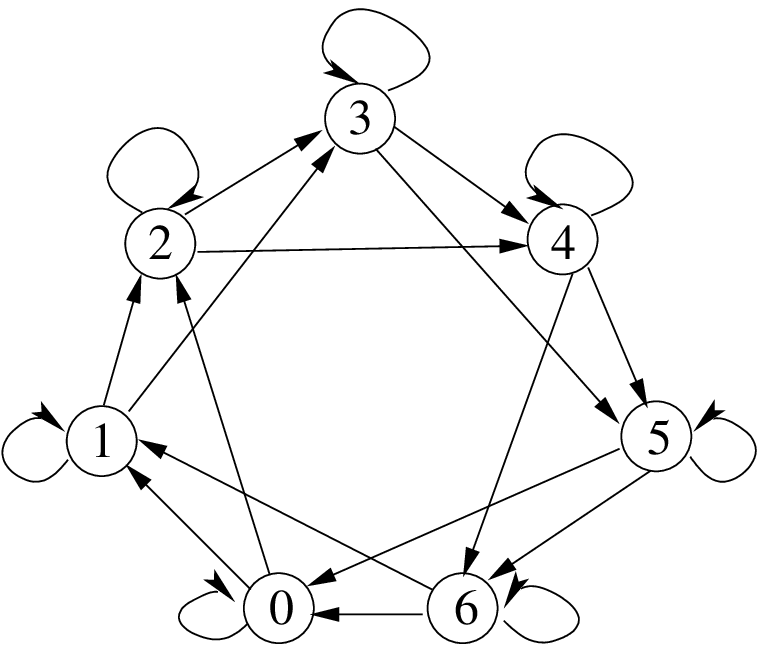, width= 2.0 in}}\\
\subfigure[$L_{6}^{0,1,2}$ and $L_{7}^{0,1,2}$]{\epsfig{file=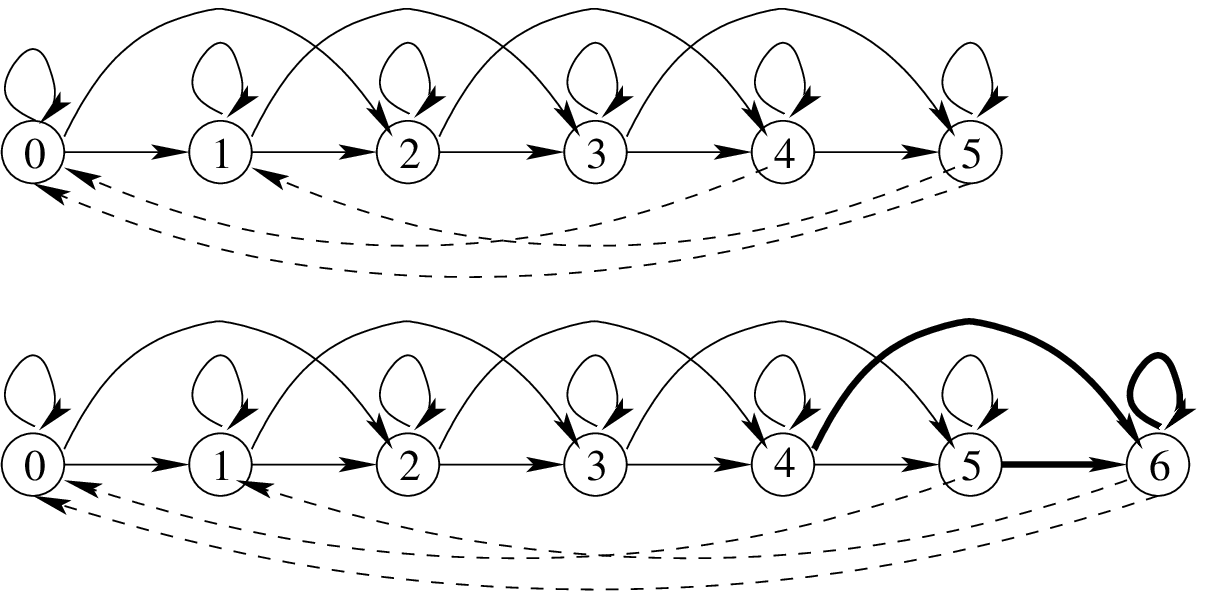, width=3.5 in}}
\end{array}
$$
\caption{$C_{n}^{0,1,2}$, a constant jump circulant and its lattice
equivalents. In (c) solid edges are $L_n;$
dashed edges are $\Hookn;$  bold
edges in $C_{7}^{0,1,2}$ are $\Nn$ for $n=7.$}
\label{fig:constant_cir}
\end{figure*}

We assume that $0 = s_1 < s_2 < \cdots < s_k$ and set $\ems = s_k.$
Figure (\ref{fig:constant_cir}) shows two circulant graphs
with constant jumps 0, 1, 2.
Note  that our assumption is
without loss of generality, as we can choose any row
of a circulant matrix to be the top one; for our assumption to be correct,
we choose a row with a '1' in its first position. Equivalently, multiplying a circulant matrix by $P_n$ or $P^{-1}_n$
doesn't change its permanent so we can normalize $S_1=0$. 
For example, $P^{-2}_n+P^{-1}_n+I$, $P^{-1}_n+I+P_n$ and $I+P_n+P^2_n$, corresponding respectively,  to graphs
$C_n^{-2,-1,0},$ 
$C_n^{-1,0,1}$
and
$C_n^{0,1,2}$,
all have
the same permanent. 

\begin{Definition}
\label{def:Lat_const}
See Figure 
\ref{fig:constant_cir}.
Let $S=\{s_1,s_2,\ldots s_k\}$, where the $s_i$ are fixed integers.
Define the  $n$-node {\em lattice graph}\footnote{The reason for calling this
a {\em lattice} graph will become visually obvious later in Definition
\ref{def:Lat_non_conts}, which generalizes this definition to the non-constant 
jump case.}
with jumps $S$ by 
$$L^{s_1,s_2,\dots,s_k}_n=(V(n),E_L(n))
$$
where 
$$E_L(n)=\Big\{(i,j) : j-i  \in S\Big\}.
$$
Now set
$$ \Hookn =  E_C(n) - E_L(n)$$
and
$$\Nn    = E_L(n+1)- E_L(n).
$$
Note that this implies
\begin{equation}
\label{eq:const_decomp}
\hspace*{.18in} 
\Lnpo = \Ln \cup \Nn
\quad
\mbox{and}
\quad
\Cn = \Ln \cup \Hookn.
\end{equation}
\end{Definition}
The simple but important observation is that, {\em when $n$ is viewed as a label
rather than as a  number},  $\Hookn$ and $\Nn$ are {\em independent} of the actual
{\em value} of $n$.

\begin{Lemma}
\label{lem:RecC_const}
\begin{eqnarray*}
\Hookn  
&=& \bigcup_{s \in S}
\left\{\, (n-j,\,s-j)  \,:\, 1 \le j \le s\right\},\\
\Nn &=& \hspace*{-.02in} \bigcup_{s \in S} \{(n-s,n)\}.
\end{eqnarray*}
Set $\ems = s_k$. Now define
\begin{eqnarray*}
L(n)&=& \{0,\ldots \ems-1\},\\
R(n)&=& \{n-\ems,\ldots,n-1\}.
\end{eqnarray*}
Then
\begin{eqnarray}
\label{eq:HN_Rec}
\Hookn &\subseteq& 
\left(R(n) \times L(n)\right) \\ 
\Nn    &\subseteq&  \left(R(n) \times \{n\}\right) \cup \left\{(n,n)\right\} 
\label{eq:Nn_Rec}
\end{eqnarray}
\end{Lemma}


{\em Important Note: In this section and the next we will always assume that $n \ge 2\ems$
since this will guarantee that $L(n) \cap R(n) = \emptyset.$ Without this
assumption some of our proofs would fail. Also note that the $\left\{(n,n)\right\}$ term in $\Nn$ appears because $0 \in S.$}

\section{A New Proof of Minc's result}
\label{sec:Minc}

Let $CC$ be a cycle-cover of $C_n$, i.e.,
$\forall v,\, \ID{T}v=\OD{T}v=1.$
  Then, from (\ref{eq:cycle_condition}), in $T=CC- \Hookn$, almost all vertices $v$ except (possibly) some of those that
have an edge of $\Hookn$  hanging off of them, have $\ID{T}v=\OD{T}v=1.$ This motivates 
\begin{Definition}
\label{def:legal_cov}
$T \subseteq E_L(n)$ is  a {\em legal cover of $L_n$} if 
\begin{itemize}
\item $\forall v \in V,\quad  \ID T v \le 1 \mbox{ and }\  \OD T v \le 1$.
\item $\forall v \in V - L(n),\quad \ID T v =1$.
\item $\forall v \in V - R(n),\quad \OD T v =1$.
\end{itemize}
\end{Definition}
Then, from (\ref{eq:const_decomp}) we have
\begin{Lemma} \label{lem:class_mot}\ \\
(a) If $T \subseteq E_C(n)$ is a cycle-cover of $C_n$, then\\
\hspace*{.4in}  
$T-\Hookn$ is a legal-cover of $L_n$.\\[0.1in]
(b) If $T \subseteq E_L(n+1)$ is a legal-cover of $\Lnpo$, then \\
\hspace*{.4in} 
$T-\Nn$ is a legal-cover of $L_n$.
\end{Lemma}

\begin{figure*}[t]
$$\begin{array}{c}
  \vspace{.2in}
\subfigure[$CC_1$ (all edges) and $T_1$ (solid edges)]{\epsfig{file=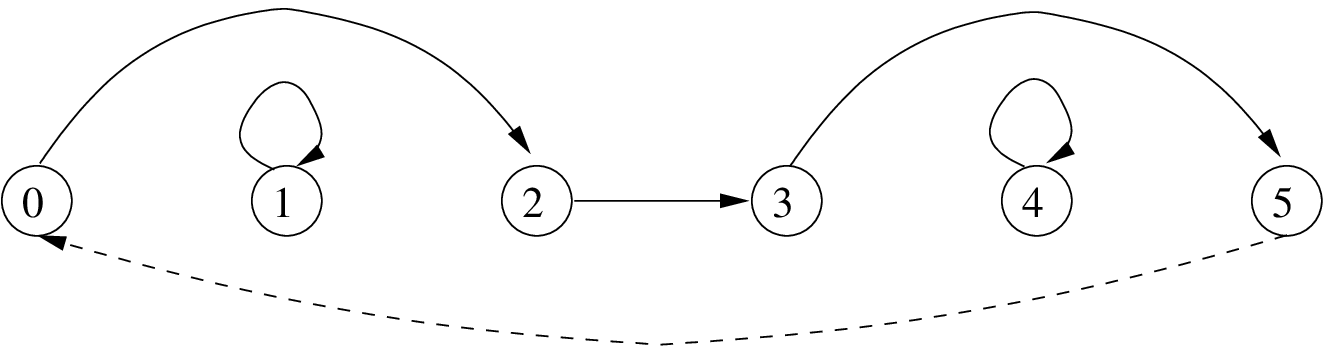, width= 2.8 in}} \\
  \vspace{.2in}
\subfigure[$CC_2$ (all edges)  and $T_2$ (solid edges)]{\epsfig{file=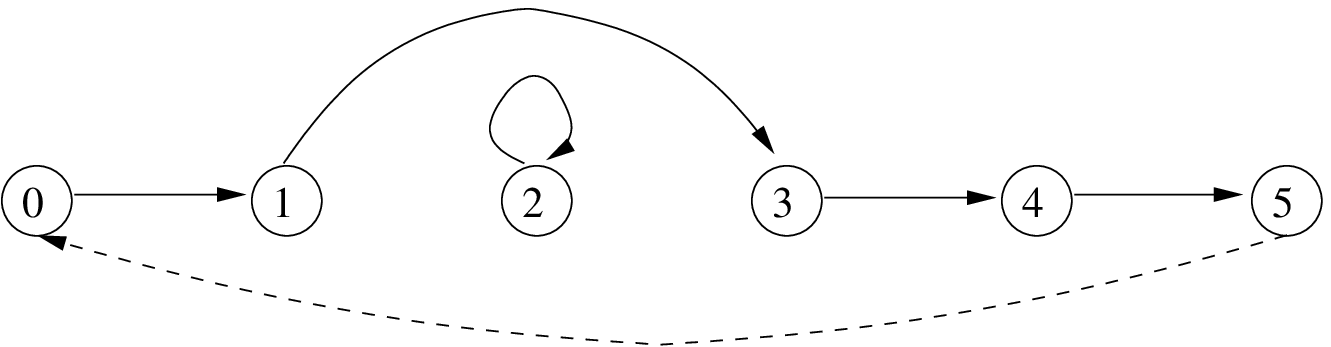, width=2.8 in}} \\
\subfigure[$CC_3$ (all edges) and $T_3$ (solid edges) ]{\epsfig{file=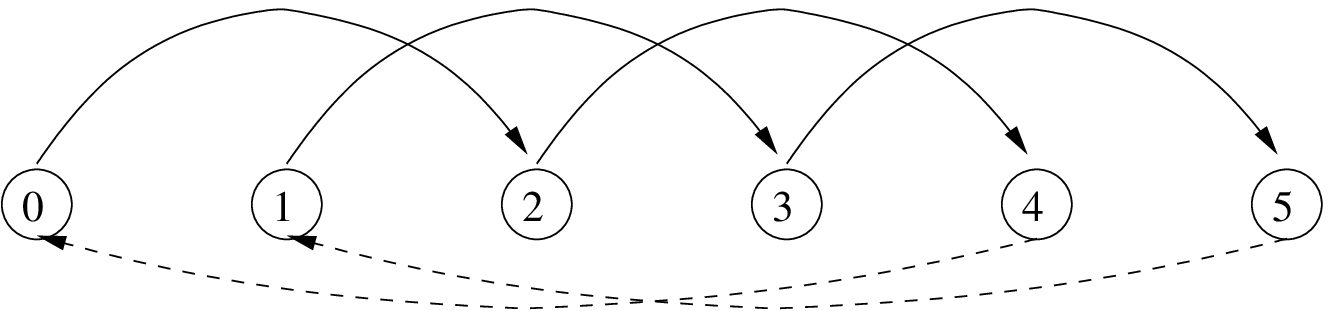, width= 2.8 in}}
\end{array}$$
\caption{All of the figures are subsets of  $C_{6}^{0,1,2}$. Solid
edges are in $L_n$; dashed edges are in $\Hookn.$ 
The solid plus dashed edges comprise three different cycle covers $CC_i,$ $i=1,2,3$ in
$C_6.$
Removing the dashed $\Hookn$ edges leaves three {\em legal} covers
$T_i,$ $i=1,2,3$, in $L_6.$
Note that $\ems=2$, $C(T_1)=C(T_2) = (\,(0,1),\,(0,1)\,)$
and $C(T_3) = (\,(0,0),\,(0,0)\,)$.
}
\label{fig:ex1}
\end{figure*}

\begin{figure*}[h]
$$\begin{array}{c}
  \vspace{.2in}
\subfigure[$T_1 \cup \{(n-2,n)\}$]{\epsfig{file=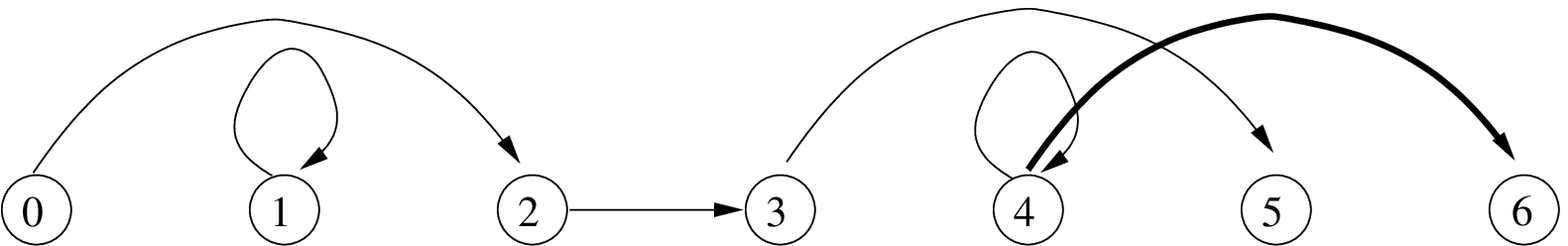, width= 3.2 in}} \\
  \vspace{.2in}
\subfigure[$T_2 \cup \{(n-2,n)\}$]{\epsfig{file=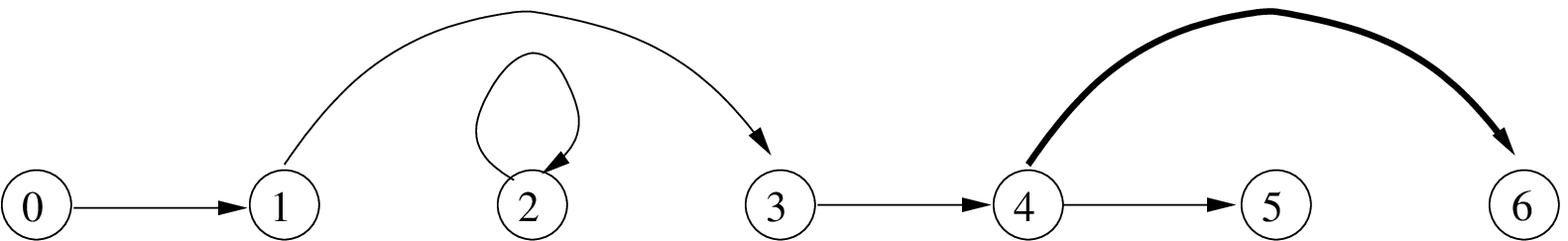, width=3.2 in}} \\
\subfigure[$T_3 \cup \{(n-2,n)\}$]{\epsfig{file=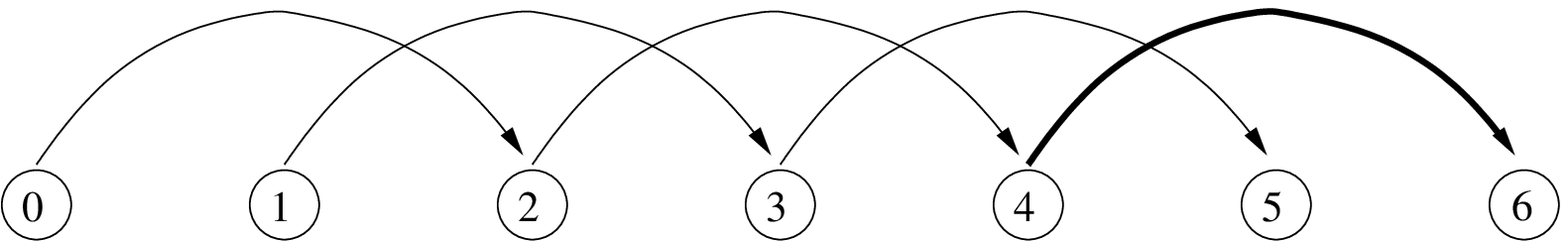, width= 3.2 in}}
\end{array}$$
\caption{$n$ was increased from $6$ to $7$ and 
 $S=\{(4,6)\} \subseteq \mbox{New}(6)$ was added to the $T_i$ of the previous figure.
Note that, in $L_7,$ $C(T_1 \cup S) = C(T_2 \cup S) = \emptyset$ since they
are no longer legal covers. 
Also, $C(T_3 \cup S) = (\,(0,0),\,(0,0)\,)$.
}
\label{fig:ex2}
\end{figure*}

 From the definition of legal covers we can classify and partition legal covers
by the appropriate in/out degrees of their  vertices in $L(n), R(n).$
\begin{Definition}
\label{def:tuple_def}
$A$ is a {\em binary $r$-tuple} if \\
\centerline{$A=(A(0),A(1),\ldots, A(r-1))$ where 
$\forall i, A(i) \in \{0,1\}$.}
\end{Definition}

\begin{Definition} (See Figure \ref{fig:ex1}).
\label{def:class_def}
Let $\calp$ be the set of $2^{2\ems}$ possible binary tuple pairs $(L, R)$ where each of 
$L, R$ are, respectively, binary $\ems$ tuples.

\vspace*{.1in}

Let $T$ be a {\em legal-cover} of $L_n.$ The {\em classification of $T$} will be 
$C(T)=   (L^T, R^T) \in \calp$
where
$$
\begin{array}{lrcl}
\forall 0 \le i < \ems,\, & L^T(i) &=& \ID T i\\[.04in]
                          & R^T(i) &=& \OD T {n-1-i}.
\end{array}
$$
If $T$ is {\em not} a legal-cover then we will write
$C(T) = \emptyset.$
Finally, set
\begin{eqnarray*}
{\cal L}(n) &=& \{T \subseteq E_L(n)\,:\, \mbox{$T$ is a legal cover of $L_n$}\}\\
{\cal L}_X(n)&=& \{T \in {\cal L}(n) \,:\, C(T) = X\}\\
T_X(n)       &=& |{\cal L}_X(n)|
\end{eqnarray*}
so $T_X(n)$ is the number of legal-covers of $L_n$ with classification $X.$
\end{Definition}
The main reason for introducing these definitions is that checking whether a legal cover $T$ of
$\Ln$ can be completed to a cycle-cover of $\Cn$ or to a legal cover in $\Lnpo$ 
doesn't depend upon  all of $T$ but {\em only upon its classification
$C(T)$}.

\begin{Lemma} \label{lem:class} See Figures \ref{fig:ex1} and \ref{fig:ex2}. \\
 Let $ X=(L^X, R^X) \in \calp$. Let $T_1$ be a legal cover in $L_{n_1}$ and $T_2$ be a legal cover of $L_{n_2}$, such that
$C(T_1)=C(T_2) =X$.

(a) Let   $S \subseteq \Hookn.$   Then, \vspace*{.1in}\\
\centerline{$T_1 \cup S$ is a cycle-cover of $C_{n_1}$}\\
\centerline{ {\bf if and only if}}\\
\centerline{ $T_2 \cup S$ is a cycle-cover of $C_{n_2}$} \\

\vspace*{.1in}

\noindent(b) 
Let  $S \subseteq \Nn.$ Then, 
$$C(T_1 \cup S) = C(T_2 \cup S).$$
That is,  either {\em both} $T_1 \cup S$ and $T_2 \cup S$ are {\em not} legal  covers
or, they are both legal covers and there is some $X'  \in \calp$ such that
$C(T_1 \cup S) = C(T_2 \cup S) = X'$
\end{Lemma}
\begin{Proof}
To prove (a) recall that $T \cup S$ is a cycle-cover of $\Ln$ if and only if,\\ 
$\forall v \in V,\, \ID {T \cup S} v = \OD {T \cup S} v =1$
or 
\begin{equation}
\forall v \in V, 
\ID S v = 1 - \ID {T} v  
\quad \mbox{and} \quad
\OD S v = 1 - \OD {T} v 
\end{equation}
 From Lemma \ref{lem:RecC_const} and the definition of
a legal cover we have that this is true if and only if
$$\begin{array}{lrcl}
  \forall i < \ems,\, & \ID S i &=& 1- L^X(i),\\[0.04in]
                        &  \OD S {n-1-i} &=& 1- R^X(i).
\end{array}
$$
and this is only dependent upon $X$ and $S$ and not upon $n$ or any other properties of $T.$

The proof of (b) is similar and omitted here.

\end{Proof}

This lemma permits us, for $X, X' \in \calp$ and $S \subseteq \Hookn$,
to abuse the notations and write {\em $(X \cup S)=X'$} to denote that,
when $C(T) = X$, $C(T \cup S)=X'$. We will sometimes also write
``{\em $X \cup S$ is a cycle cover}'' to denote that $T \cup S$ is a 
cycle cover.


%


\begin{Definition} 
\label{def:alph_bet}
For $X,X' \in \cal P,$ $S \subseteq \Hookn$ and $S' \subseteq \Nn$ set
$$
\beta_{X,S} = \left\{
\begin{array}{ll}
1 & \mbox{ if $X \cup S$ is  a cycle cover}\\
0 & \mbox{ otherwise}
\end{array}
\right.
$$
and 
$$
\alpha_{X,X',S'} = \left\{
\begin{array}{ll}
1 & \mbox{ if $C(X' \cup S')=X$}\\
0 & \mbox{ otherwise}
\end{array}
\right.
.
$$
Now set
$$
\beta_X = \sum_{S \subseteq \Hookn} \beta_{X,S}
$$
and
\begin{equation}
\label{eq:alpha_def}
\alpha_{X,X'} = \sum_{S' \subseteq \Nn} \alpha_{X,X',S'}.
\end{equation}
\end{Definition}
Note that $\beta_X$ and $\alpha_{X,X'}$ are constants that can be mechanically calculated.
In fact $\alpha_{X,X'}$ is much simpler to calculate than  it might initially appear seem since
\begin{Lemma}
\label{lem:alpha_calc}
If $\alpha_{X,X',S'} =1$, then  $|S'| =1$.
\end{Lemma}
\begin{Proof}
In order for $X'\cup S$ to be a legal cover $S$ must include at least one edge that
points to vertex $n$, so $|S| \ge 1.$
From (\ref{eq:Nn_Rec}),   {\em all} edges in $\Nn$ point to $n$.  If $|S'|>1$,
then 
$\ID{X'\cup S} n  = |S'|>1$  and $X'\cup S$  wouldn't be a legal cover.
\end{Proof}

Thus (\ref{eq:alpha_def}) can be calculated by summing over $|\Nn| = k$ values,
instead of $2^{k}$ values.

Lemmas \ref{lem:class_mot}  and \ref{lem:class} immediately imply  our first technical result.
\begin{Lemma}
\label{lem:main_const}
$$T(n) = \sum_{X \in \calp} \beta_X T_X(n)$$
and
$$
T_X(n+1) = \sum_{X' \in \calp} \alpha_{X,X'} T_X(n).
$$
\end{Lemma}

Let $m = |\calp| = 2^{2 \ems}.$ Take any arbitrary ordering of $\calp$
and define the $1 \times m$ constant vector $\beta= (\beta_X)_{X \in \calp}$ and 
$m \times m$ constant matrix $A = (\alpha_{X,X'})_{X,X' \in \calp}.$  Finally, set
$\bar T(n)= \mbox{col}(T_X(n))_{X \in \calp}$ to be a $m \times 1$ column vector.
Then, Lemma \ref{lem:main_const} is exactly
\begin{eqnarray}
\nonumber
\forall n \ge 2\ems,& &
T(n) = \beta\,   \bar T (n)
\quad \mbox{and} \quad 
\bar T(n+1) =A \, \bar T(n)
\end{eqnarray}
which is equation (\ref{eq:trans_intro}). As mentioned in the introduction, this
immediately implies that $T(n)$ satisfies a fixed-degree constant coefficient recurrence
relation where the degree of the recurrence is at most the degree of any polynomial
$P(x)$ such that $P(A)=0.$ By the Cayley-Hamilton theorem, $Q(A)= 0,$
$Q(x)$ is the  degree $m=2^{2 \ems}$
characteristic polynomial $Q(x) = \mbox{det}(IX-A)$.

We will now see that it is possible to reduce this degree from $2^{2 \ems}$ down to below $2^{\ems}.$  We will do this by showing that,  given  appropriate orderings of 
the classifications,  $A = (\alpha_{X,X'})$ will have a very special block diagonal
format.  In what follows, please refer to the worked example in Appendix 
\ref{App:worked example} for illustration.

\begin{Definition}
\label{def:consistent}
A linear ordering on the classifications  $\calp$ will be called {\em consistent} if it is the lexicographic concatenation
of linear orderings on its left and right components. 

More specifically, linear ordering ``$<$'' is consistent on $\calp$ if there exist linear orderings
``$\Lleq$'' and ``$\Rleq$'' such that if $X_1=(L_1,R_1)$ and $X_2=(L_2,R_2)$ we have
$X_1< X_2$ if and only if one of the following is true
$$ L_1 \,\Llt\, L_2
\quad\quad
\mbox{or}
\quad\quad 
L_1\,\Leq\,  L_2 \ \mbox{ and  } \ R_1 \,\Rlt\, R_2
$$
\end{Definition}

Note that in the above definition it is not necessary for the ordering on the left component to be
the same as the ordering on the right one (we will use this  fact later in Lemma \ref{lem:degree_bound}).

\begin{Lemma}
\label{lem:high_deg_bound}
Let $A = (\alpha_{X,X'}).$ If $X \in \calp$ is ordered consistently, then 
\begin{equation}
\label{eq:defAbar}
A  = 
\left(
   \begin{array}{cccc}
      \bar A & 0 & \cdots & 0 \\
      0 & \bar A & 0 & 0 \\
      \vdots & 0 & \ddots & \vdots \\
      0 & 0 & \cdots & \bar A
   \end{array}
\right)
\end{equation}
where $\bar{A}$ is some $2^\ems \times 2^\ems$ matrix.  That is, 
 $A=\mbox{diag}(\bar A,\bar A,\ldots,\bar A)$ where $A$ contains $2^{\ems}$ copies of $\bar A$ on its diagonal.

\end{Lemma}
\begin{Proof}
Suppose $X=(L^{X}, R^{X})$ and $X'=(L^{X'}, R^{X'})$.  

Recall 
that $\alpha_{X,X'} = \sum_{S \subseteq \Nn} \alpha_{X,X',S}$ where
$\alpha_{X,X',S}=1$ if and only if $C(X' \cup S)=X$, and is otherwise $0.$

\hspace*{.1in}

Let $L$ denote any  binary  $\ems$-tuple. Partition
$\calp$ up into 
$2^{\ems}$ sets of size $2^{\ems},$ 
$\calp_{L}  = \{ X \in \calp \,:\, L^{X} = L\}$.

\hspace*{.1in}

Note that, from Lemma \ref{lem:RecC_const}, if  $S \subseteq \Nn$, none of $S$'s  edges have endpoints in $L(n)$.
Intuitively,  this is because edges in $\Nn$ only connect vertices near the
{\em right} side of the lattice and do not touch any vertices on the {\em left} side of
the lattice.

Thus,  if $\alpha_{X,X',S}=1$, then $L^{X} = L^{X'}$.
In particular this means that if $\alpha_{X,X',S}=1$ then $X,X'$ are both in
the same partition set $\calp_{L}$.


\hspace*{.1in}

Now suppose that  $\alpha_{X,X',S}=1.$
Let  $\bar L$  be {\em any} other binary $\ems$-tuple and set
\begin{equation}
\label{eq:same}
\bar X=(\bar L, \, R^{X})
\quad\mbox{and}\quad
\bar X'=(\bar L, \, R^{X'}).
\end{equation}
Then, 
again using the fact that none of the endpoints of $S$ are in $L(n)$
we have  that
$C(X' \cup S)=X$ if and only if $C(\bar X' \cup S)=\bar X$ so $\alpha_{X,X'} = \alpha_{\bar X,\bar X'}.$

\hspace*{.1in}

When constructing matrix $A = (\alpha_{X,X'})_{X,X'\in \calp}$ we previously allowed any
arbitrary ordering of $\calp.$  Ordering 
 the $X \in \calp$ consistently   groups all of the
$X$ in a particular $\calp_{L}$ consecutively.  The observations above imply that
$A$ is partitioned into $2^{\ems} \times 2^{\ems}$ blocks where each block
is of size $2^{\ems} \times 2^{\ems}$.  The non-diagonal blocks correspond to
$\alpha_{X,X'}$ where $X,X'$ are in different partitions so all of the non-diagonal blocks
are $0.$  On the other hand, the fact that $\alpha_{X,X'} = \alpha_{\bar X,\bar X'}$ for the
$\bar X, \bar X'$ defined in (\ref{eq:same}) and the consistency of the ordering of the $X$
tells us that
all the diagonal blocks are copies of each other, i.e., 
we have proven (\ref{eq:defAbar}).

\end{Proof}

\begin{Cor}
There is a degree $2^\ems$ polynomial $P(x)$ such that $P(A)=0$.
\end{Cor}
\begin{Proof} 
From the previous lemma,
any polynomial $P(x)$ that annihilates $\bar A$ also annihilates $A.$
Since $\bar A$ is a  $2^{\ems} \times 2^{\ems}$ matrix,
the  Cayley-Hamilton theorem says that
the characteristic polynomial $\bar P(x)$ of $\bar A$, which is of degree $2^{\ems}$,
annihilates $\bar A$.
\end{Proof}

The original Minc result \cite{Minc78,Minc87b}) gave an order of
$2^{\ems}-1$.  We can derive this through a slightly more sophisticated decomposition
of $\bar{A}.$

\begin{Lemma}
\label{lem:degree_bound}
Let $A = (\alpha_{X,X'}).$ 
Then there is a degree $2^{\ems}-1$ polynomial
$P(x)$ such that $P(A)=0.$
\end{Lemma}
\begin{Proof}
If $\alpha_{X,X',S}=1$, then we have just seen that $L^X=L^{X'}$. 
Consider $R^X$, $R^{X'}$ and $S$.

First recall from Lemma \ref{lem:alpha_calc} that if $\alpha_{X,X',S}=1$
then 
$S$ contains exactly  one edge.

We claim that if $\alpha_{X,X',S}=1$ then  
 $X$  and $X'$ must contain exactly the same
number of `$0$'s. Note that since $L^X = L^{X'}$ we only need to show that
$R^X$ and $R^{X'}$  have the same number of `$0$'s.

There are actually two cases.
The first case is that $S=\{(n-s,n)\}.$
In this case we are throwing away one vertex ($n-s$) which (because of the legality of $X'$)
had outdegree zero
and adding  a new vertex $n$ which also has outdegree $0.$
So, the number of `$0$'s in $X \cup S$ is the same as the
number of `$0$'s in $X'.$

The second case is that  $S=\{(n-i,n)\}$ where $i < s.$.
Since $X' \cup S$ is legal, vertex $n-s$ must have already had outdegree one so 
throwing it away doesn't change the number of `$0$'s. 
Adding the new vertex $n$ with outdegree '$0$' increases the number of `$0$'s by one.
Adding edge $(n-i,n)$ changes the outdegree of vertex $n-i$ to one,
decreasing the number  of `$0$'s   by one.

So, the number of `$0$'s in $X' \cup S$ is again the same as the
number of `$0$'s in $X.$

Recall that we have that $L^X=L^{X'}.$
This suggests that we can re-order the entries of $\bar A$
so that all $R^X$ with the same number of 0's are grouped together (maintaining 
the fact that the ordering is consistent).
Since there are ${\ems} \choose i$ $\ems$-tuples containing $i$ `$0$'s,
$\bar A$ will become a block diagonal matrix
of $s+1$ blocks with block $i$ having size ${\ems} \choose i$.
That is
\begin{equation}
\label{eq:defB}
 \bar{A} = 
\left(
   \begin{array}{ccccc}
       B_0 & 0 & \cdots & 0 & 0\\
      0 & B_1 & \cdots & 0 & 0\\
      \vdots & \ddots & \ddots & \ddots & \vdots \\
      0 & 0 & \cdots & B_{n-1} & 0\\
      0 & 0 & \cdots & 0 & B_n
   \end{array}
\right)
\end{equation}
where $B_i$ is  a  
${{\ems} \choose i} \times {{\ems} \choose i}$
matrix..

Let $P_i(x)$ be the characteristic polynomial of $B_i.$  This has degree $\le {{\ems \choose i}}$.

Note that  $B_0$ and $B_n$ are both $1 \times 1$   matrices. 
By construction, 
 $B_0= B_n =(1)$ so 
$P_0(x) = P_n(x) = 1-x$, i.e., their characteristic polynomial is the same..

Because of the block diagonal form of $B,$  $P(x) = \prod_{i=0}^{n-1} P_i(x)$ annihilates
$B.$  This polynomial has degree
$\le \sum_{i=0}^{n-1} {{\ems \choose i}} = 2^{\ems} -1$, proving the lemma.
\end{Proof}

Lemma  \ref{lem:main_const} tells us that 
(\ref{eq:trans_intro}) holds while Lemmas  \ref{lem:high_deg_bound} and \ref{lem:degree_bound}  tell us that
matrix $A$ is annihilated by polynomial $P(x)$ of degree $2^{\ems}-1$.  Combining
them  gives that $T(n)$ satisfies a 
degree-$(2^{\ems}-1)$  constant coefficient recurrence relation.

\subsection{Deriving the Recurrence Relation}


We have just seen that $T(n)$ {\em satisfies} a 
degree-$(2^{\ems}-1)$  constant coefficient recurrence relation. 
To actually {\em derive} the recurrence relation we must construct 
\begin{itemize}
\item {(i)} a polynomial 
$Q(n)$
that annihilates $A = \{\alpha_{X,X'}\}$ and 
\item (ii) the initial conditions
$T(n),$  $n = 2 \ems,\,  2 \ems + 1,\, \ldots, \, 2 \ems + 2^{\ems} -2.$
\end{itemize}


To construct $Q(x)$, note from Lemma  \ref{lem:high_deg_bound} 
that it suffices to calculate the characteristic polynomials $Q(x)$
of matrix $\bar{A}.$
We must therefore first calculate the 
$2^{2\ems}$  $\alpha_{X,X'}$ entries of  $\bar{A}.$  


Recall that 
$
\alpha_{X,X'} = \sum_{S \subseteq \Nn} \alpha_{X,X',S}
$
and, as noted in the  proof of Lemma \ref{lem:degree_bound}, we know that
if $\alpha_{X,X',S'} =1$, then $S$ contains at most one edge.  Since
$\Nn$ contains $\ems$ edges, we can,  with
the appropriate data structures,  calculate $\alpha_{X,X'}$ in $O(\ems^2)$ time.
We can therefore calculate all the non-zero entries in the $\bar{A}$  in
$O(\ems^2 2^{2\ems})$ time.
Finally, we can calculate  $Q(x)$ in 
$O(2^{3\ems})$ 
 time, since it takes $O(n^3)$ time to compute the
characteristic polynomial of an $n \times n$ matrix \cite{Householder75}.
Thus,  we can calculate $Q(x)$ in $O(2^{3\ems})$ time.

To derive (ii), the initial conditions
$T(n),$ 
$n = 2 \ems,\,  2 \ems + 1,\, \ldots,\, 2 \ems + 2^{\ems} -2.$
suppose first that we already knew
$\bar T(2 \ems)$ and $\beta$.  Since $\bar T(n+1) = A \bar T(n)$ we can use
the block structure from (\ref{eq:defAbar}) to calculate $\bar T(n+1)$ from $\bar T(n)$ in 
$O(2^{3 \bar s})$ time.  It then takes only another $2^{2 \bar s}$ time to calculate
$T(n+1) = \beta \bar T(n+1).$  So,  we can calculate all of the values 
$T(n),$  $n = 2^{\ems}+1, 2^{\ems}+1,\ldots, 2\ems + 2^{\ems} -2$ in $O(2^{4 \ems})$ time,
improving upon the doubly exponential procedure implied by
Minc's original result.


It still remains to calculate $T_X(2 \ems)$ and $\beta_X$ for all classifications $X$.

Let $X=(L,R).$ We want to calculate the number of legal covers in $L_{2\ems}$ with
classification $X.$ In a legal cover the number of '$0$'s in $L$ must be  equal to the
number of '$0$'s in $R$.  Let $a_1,a_2,\ldots,a_i$ be the indices such that $L(i)=0$  and
$b_1,b_2,\ldots,b_i$ be the indices such that $R(i)=0.$ Define the set of $i$ edges
$A = \bigcup_{j=1}^i \{(2\ems - 1 - b_j, a_j)\}$.   Now define a new graph $G_X$ as follows:
(a) start with the lattice graph $L_{2 \ems};$ (b) remove all edges
{\em entering} vertices $a_j$, $j=1,2,\ldots,i;$ 
(c) remove all edges {\em leaving} vertices $2 \ems -1 - b_j$,
$j=1,2,\ldots,i;$  (d) add the $i$ edges in $A.$  Then it is not
difficult to see that  
$T$ is a legal cover in $L_{2 \ems}$  if and only if $T \cup A$ is a cycle cover of $G_X.$  Since
every cycle cover of $G_X$ {\em must} contain all edges in $A$ there is a one-one
correspondence
between cycle covers in $G_X$ and legal covers in $L_{2 \ems}$ with classification $X.$
We can therefore calculate $T_X(2 \ems)$ by calculating the permanent of the adjacency matrix of
$G_X$ which can be done in $O(s 2^{2 \ems})$ time using Ryser's algorithm.  Calculating
all entries in $\bar T(2 \ems)$ then takes $O(s 2^{4 \ems})$ time.

Finally,  we must calculate all the $\beta_X$.  
Let $X=(L,R).$  If $\beta_X \not=0$ then the number of '$0$'s in $L$ must be  equal to the
number of '$0$'s in $R$.  As above, let $a_1,a_2,\ldots,a_i$ be the indices such that $L(i)=0$  and
$b_1,b_2,\ldots,b_i$ be the indices such that $R(i)=0.$  Now construct the $i \times i$  bipartite 
graph $B$ as follows:\\
\centerline{  Edge $(j,k) \in B$ if and only if $(n-1 - b_j,a_k) \in
  \Hook(2 \ems).$}\\
It is not difficult to see that $\beta_x$ is exactly the number of complete matchings in 
$B$.  We can therefore calculate $\beta_X$ by evaluating the permanent of the adjacency matrix
of $B$.  This can be done in $O(\ems 2^{\ems})$ time per entry and thus in $O(\ems 2^{3\ems})$
time in total.


Combining everything,  we see that we can construct the recurrence relation and initial
conditions using $O(s 2^{4 \ems})$ time.


\section{Non-constant Jump Circulant Graphs}
\label{sec:non_constant}
We now extend the above definitions and lemmas to the case of non-constant
circulants $\Cn = C_{pn+s}^{p_1n+s_1,\,p_2n+s_2,\,\cdots,\,p_kn+s_k}$ 
as introduced in Definition \ref{def:circ_non_const}.
Note that if $s = \alpha  p + \beta$ for some arbitrary integer $\alpha$ and integer $\beta \ge 0,$ 
we can rewrite $\Cn$ as 
$$ C_{p(n+\alpha) +\beta}^{p_1(n+\alpha)+s'_1,\,p_2(n+\alpha)+s'_2,\,\cdots,\,p_k(n+\alpha)+s'_k}$$
where $\forall i, s'_i = s_i - \alpha p$.  Thus,  we may assume $0 \le s < p.$


Note that, using  a similar argument
to that in the previous section preceding Definition
\ref{def:Lat_const}, we may and do, without loss of generality,
assume $\forall i, \, s_i \geq  s.$

\begin{figure*}[p]
$$
\subfigure[$C_{4n+1}^{1,n+2,2n+1}$ for $n=4$]{\epsfig{file=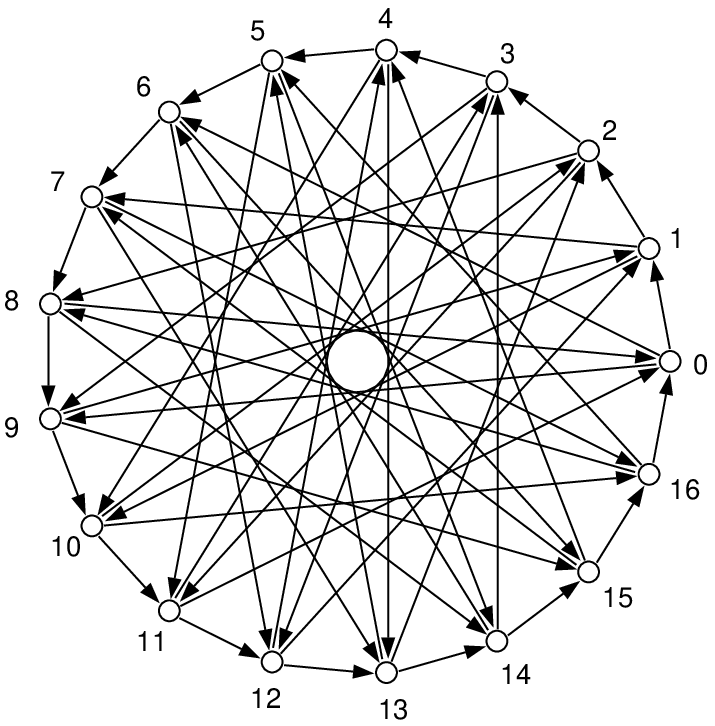, width=2.4in}}
\hspace*{.4in}
\subfigure[$C_{4n+1}^{1,n+2,2n+1}$ for $n=4$]{\epsfig{file=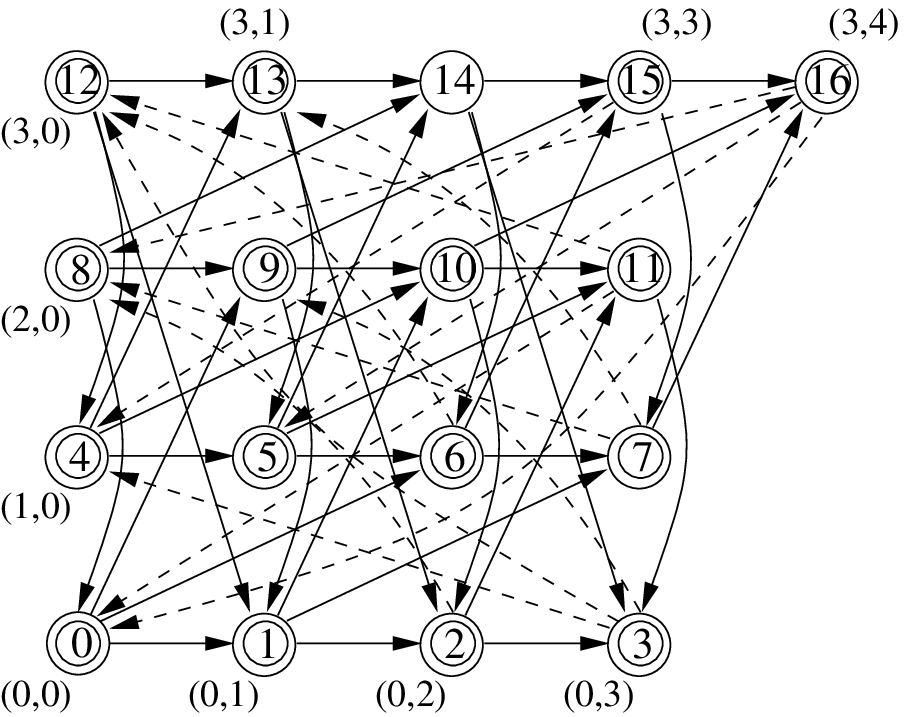, width=3.6 in}}
$$
$$
\hspace*{-.6in}
\subfigure[$C_{4n+1}^{1,n+2,2n+1}$ for $n=4$]{\epsfig{file=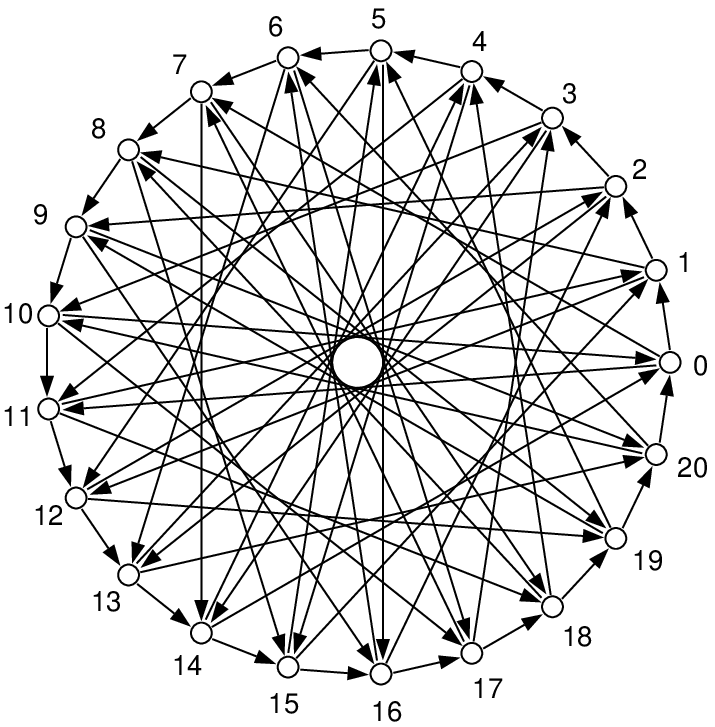, width=2.7in}}
\hspace*{.4in}
\subfigure[$C_{4n+1}^{1,n+2,2n+1}$ for $n=5$]{\epsfig{file=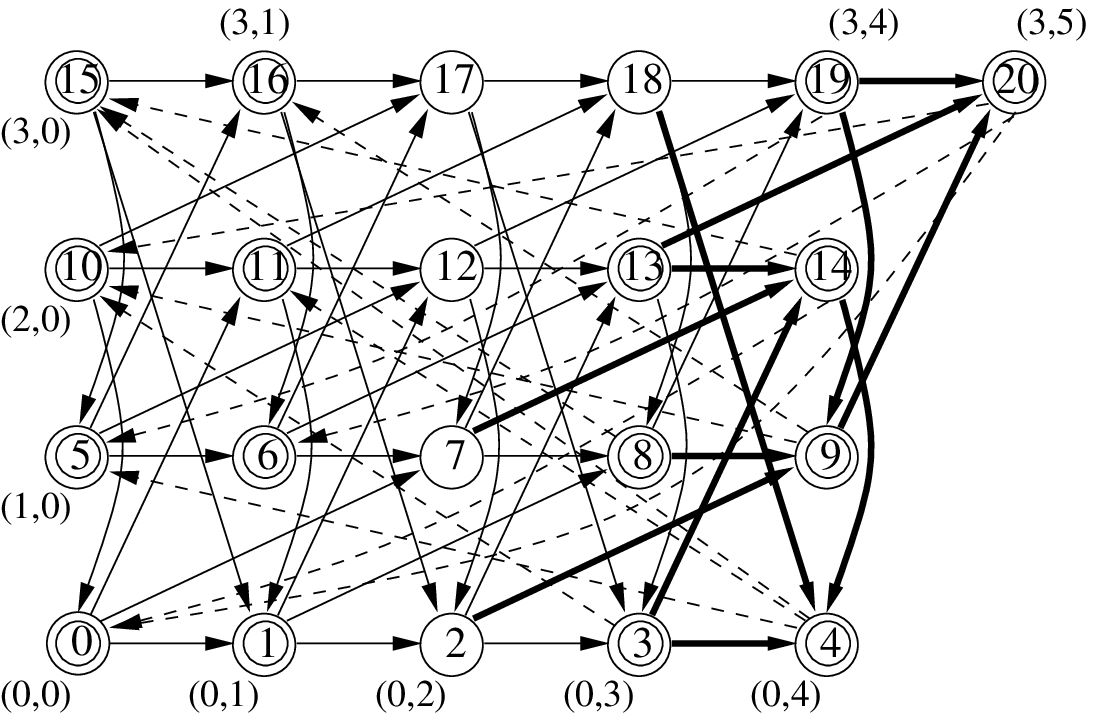, width=4.2 in}}
$$
\caption{
The graphs $C_{4n+1}^{1,n+2,2n+1}$.  (a) and (b) are two representations of  $n=4$; (c) and (d) 
are $n=5$.  (a) and (c) are drawn in traditional circulant format;  (b) and (d) in
lattice graph format. 
The bold edges in (d) are $\Nn$. In the lattice graph format  the dashed edges 
are $\Hookn$ and the 
double-circled nodes denote  $L(n)$ on the left and $R(n)$ on the right.
In (b), $L(n)$ and $R(n)$ actually abut each other; 
$L(n)$ are the double-circled nodes in the two leftmost columns;  $R(n)$
are the double-circled nodes in the three rightmost columns.
As discussed in the text,  all edges in $\Hookn$ go from $R(n)$
to $L(n).$
}
\label{fig:non_const2}
\end{figure*}

Analyzing non-constant jump circulants    will
require a change in the way that we visualize the nodes of $C_n;$  until, 
now, as in Figure \ref{fig:const}(c),  we 
visualized them as points on a line with the edges in $\Hookn$ connecting the left and
right endpoints of the line.
In the non-constant jump case it will be convenient to visualize them
as points on a bounded-height {\em lattice}, where $\Hookn$ connects the left and
right boundaries of the lattice. We start by introducing a new graph:


\begin{Definition}
\label{def: Circ Lattice}
See Figures   \ref{fig:non_const} and \ref{fig:non_const2}.
Let $p,s$,    $p_1,p_2,\ldots,p_k$  and $s_1,s_2,\ldots,s_k$
be given  non-negative integral constants such that $\forall i,$ $0
\le p_i<p$ and $s_i \geq s\ge 0.$ 
Set $S = \{p_1n+s_1,p_2n+s_2,\cdots,p_kn+s_k\}.$
For $u,v$ and integer $n$, set
$f(n;u,v)= un+v.$
Define
$${\widehat C}_n = \left({\widehat V}_C(n),\, {\widehat E}_C(n)\right)$$
where 
$$
{\widehat V}(n) =
\left\{
(u,v) \, \bigg|\, 
\begin{array}{c}
0 \le u \le p-2\\
0 \le v \le n-1
\end{array}
\right\}
\ \bigcup\ 
\left\{\, (p-1,v) \,:\, 0\leq v\leq n+s-1\right\}
$$

and
$$ {\widehat E}_C(n) = 
\left\{
((u_1,v_1), (u_2,v_2)) \, \bigg|\, 
\begin{array}{l}
(u_1,v_1), (u_2,v_2) \in \widehat V_C(n) \mbox{ and }\\
\quad \exists i \mbox{ such that}\\
 \quad f(n;u_2,v_2)-f(n;u_1,v_1) = p_i n + s_i \bmod{(pn+s)}
\end{array}
\right\}
$$
%
 %


\end{Definition}


If $s=0$ we see that  ${\widehat C}_n$ is simply a rectangular lattice with
other regularly placed edges
as in Figure \ref{fig:non_const}.
If $s >0$ then ${\widehat C}_n$
is a rectangular lattice with extra vertices extending out from its top row as
in Figure \ref{fig:non_const2}.  These extra vertices, which disturb the regularity
of the lattice,  are what will complicate our analysis.

Directly from the definition we see ${\widehat C}_n$ is {\em isomorphic} to
$\Cn = C_{pn+s}^{p_1n+s_1,p_2n+s_2,\cdots,p_kn+s_k}$.  In particular, cycle-covers of
${\widehat C}_n$ are in 1-1 correspondence with cycle covers of $\Cn$ so we can
restrict ourselves to counting cycle covers of ${\widehat C}_n$.  We now introduce the generalization
of Definition \ref{def:Lat_const}.


\begin{Definition}
\label{def:Lat_non_conts}
Let $p,s$,    $p_1,p_2,\ldots,p_k$  and $s_1,s_2,\ldots,s_k$ and $S,$ $f$ be as in Definition
\ref {def: Circ Lattice}.
Define the  $pn+s$-node {\em lattice   graph}
with jumps $S$
$$\Ln=\left({\widehat V}(n),{\widehat E}_L(n)\right)$$
%
where
$$ {\widehat E}_L(n) = 
\left\{
((u_1,v_1), (u_2,v_2)) \, \bigg|\, 
\begin{array}{l}
(u_1,v_1), (u_2,v_2) \in \widehat V_C(n) \mbox{ and }\\
\quad \exists i \mbox{ such that}\\
 \quad (a)\,  f(n;u_2,v_2)-f(n;u_1,v_1) = p_i n + s_i \bmod{(pn+s)}\\
 \quad \mbox{\bf and }\\
 \quad  (b) \, u_2-u_1=p_i \bmod{p}
\end{array}
\right\}
$$
%
Now set
$$ \Hookn =  {\widehat E}_C(n) - {\widehat E}_L(n)$$
and 
$$\Nn    = {\widehat E}_L(n+1)- {\widehat E}_L(n).
$$
Note that this implies 
\begin{eqnarray}
\label{eq:non_const_decomp}
\nonumber \Lnpo &=& \Ln \cup \Nn\\
&\mbox{and}&\\
{\widehat C}_n &=& \Ln \cup \Hookn.\nonumber
\end{eqnarray}
\end{Definition}

We need  the following intuitive lemma that was used implicitly in the constant
jump case (but was so obvious there that it was not explicitly mentioned). 
The proof is straightforward but tedious and has therefore been moved to Appendix 
\ref{sec:appB}.

\begin{Lemma}
\label{lem:L_n}
Let $(u_1,v_1), (u_2,v_2) \in \widehat{V}(n),$ and $e =((u_1,v_1), (u_2,v_2)),$

$$e \in \widehat{E}_L(n) \Leftrightarrow e \in \widehat{E}_L(n+1).$$
\end{Lemma}

In the constant jump case, we were able to define $L(n)$, $R(n)$  such that
all edges in $\Hookn$ went from $R(n)$ to $L(n)$. 
It turns out that  this property remains  in the  non-constant jump case as well. 
However, as will be
seen from the internals of the proof of Lemma~\ref{lem:RecC_non_const}, 
this property is a result of our assumption that  $s_i \geq s$ for all $i$'s.
If  this assumption did not hold,  then  some $\Hookn$ edges might go  from
$L(n)$ to $R(n)$.


It is now straightforward to derive  an analogue to Lemma \ref{lem:RecC_const} showing that
$\Hookn$ and $\Nn$ are {\em independent} of the actual
{\em value} of $n$.  Before doing so we  will  need one more definition:
\begin{Definition}
$$NV(n) = V_L(n+1)-V_L(n).$$
\end{Definition}
$NV(n)$ will be the {\em new} vertices in $V_L(n+1).$  Note that we did not
explicitly define this for {\em fixed-jump} circulant graphs since in the fixed-jump case  $NV(n)=V_L(n+1)-V_L(n) = \{n\}$, i.e., there was
only the one new vertex at each step.


\begin{Lemma}
\label{lem:RecC_non_const}
Set $\ems = s_k$, and define 
{\small
\begin{eqnarray*}
L(n) & = & \left\{ (u,v): 0 \leq u \leq p-1 \mbox{ \rm{and} } 0 \leq v \leq \ems-1 \right\} \\
R(n) & = & \left\{ (u,v): 0 \leq u \leq p-2 \mbox{ \rm{and} } n-\ems \leq v \leq n-1 \right\} \\
     &   & \quad  \cup\, \left\{ (p-1,v): n+s-\ems \leq v \leq n+s-1 \right\} 
\end{eqnarray*}
}
Then
\begin{eqnarray*}
\nonumber 
\Hookn &\subseteq& R(n) \times L(n) \\
\Nn    &\subseteq& \left( R(n) \times NV(n) \right) \cup \left( NV(n) \times NV(n) \right).
\end{eqnarray*}
\end{Lemma}
The proof is straightforward but tedious and has therefore also  been moved to Appendix \ref{sec:appB}.
Figures \ref{fig:non_const} and \ref{fig:non_const2} illustrate the lemma.

In Section \ref{sec:Minc} we described 
how to calculate the number of cycle-covers in constant-jump circulant graphs.
Reviewing the proof, everything there followed directly as a consequence from the recursive
decomposition of circulant graphs in (\ref{eq:const_decomp}) combined with the structural
properties of the decomposition given in Lemma \ref{lem:RecC_const}.  But, as we have just seen,  non-constant jump circulants and their decompositions have exactly the {\em same} structural properties,
given in (\ref{eq:non_const_decomp}) and Lemma \ref{lem:RecC_non_const}.  Therefore,
the entire proof  developed in Section \ref{sec:Minc}  can be rewritten to work
for non-constant jump
circulants. The equivalent definitions and lemmas  needed in the
non-constant jump case are stated below.

\begin{Definition}
$T \subseteq {\widehat E}_L(n)$ is  a {\em legal cover of $L_n$} if 
\begin{itemize}
\item $\forall v \in V,\quad  \ID T v \le 1 \mbox{ and }\  \OD T v \le 1$.
\item $\forall v \in V - L(n),\quad  \ID T v =1$.
\item $\forall v \in V - R(n),\quad  \OD T v =1$.

\end{itemize}
\end{Definition}

\begin{Lemma} \ \\
(a) If $T \subseteq {\widehat E}_C(n)$ is a cycle-cover of $C_n$, then\\
\hspace*{.4in}    $T-\Hookn$ is a legal-cover of $L_n$.\\
(b) If $T \subseteq {\widehat E}_L(n+1)$ is a legal-cover of $\Lnpo$, then\\
\hspace*{.4in}    $T-\Nn$ is a legal-cover of $L_n$.
\end{Lemma}

The only major rewriting is required  in the analogue to Definition
\ref{def:class_def}.  The more complicated structure of the
lattice graph in the non-constant jump case requires a more complicated function to
map the indices of the $L(n)$ and $R(n)$ nodes.
\begin{Definition}
\label{def:cycle_to_legal_non_const}
$A$ is a {\em binary $r$-tuple} if \\$A=(A(0),A(1),\ldots, A(r-1))$ where 
$\forall i, A(i) \in \{0,1\}$.
\vspace*{.1in}
Let $\calp$ be the set of $2^{2p\ems}$ tuples $(L, R)$ where 
$L, R$ are two binary $p\ems$ tuples.
\vspace*{.1in}
Let $T$ be a {\em legal-cover} of $L_n.$ The {\em classification of $T$} will be 
$C(T)=   (L^T, R^T) \in \calp$
where
$$
\begin{array}{lrcl}
\forall 0 \le i < p\ems,\, & L^T(i) &=& \ID T {g^L(i)} \\[.04in]
                           & R^T(i) &=& \OD T {g^R(i)} \\[.04in]
\end{array}
$$
where
$$
\begin{array}{rcl}
\vspace{.1in}
g^L(i) & = & (\lfloor i/\ems \rfloor, i \mbox{ \rm mod } \ems) \\
\vspace{.1in}
g^R(i) & = & \left\{ \begin{array}{ll} (\lfloor i/\ems \rfloor, n-1-(i \mbox{
      \rm mod } \ems)) & \lfloor i/\ems \rfloor < p-1 \\ (\lfloor i/\ems
    \rfloor, n+s-1-(i \mbox{ \rm mod } \ems)) & \mbox{otherwise} \end{array}
\right. 
\end{array}
$$
{\em \small Note: $g^L$ and $g^R$ are simply mappings of the indices of
the $L^T(i)$ and $R^T(i)$ tuples to the nodes in $L(n)$ and $R(n)$.}\\
If $T$ is {\em not} a legal-cover then we will use the convention that
$C(T) = \emptyset.$
Finally, set
\begin{eqnarray*}
{\cal L}(n) &=& \{T \subseteq E_L(n)\,:\, \mbox{$T$ is a legal cover of $L_n$}\}\\
{\cal L}_X(n)&=& \{T \in {\cal L}(n) \,:\, C(T) = X\}\\
T_X(n)       &=& |{\cal L}_X(n)|
\end{eqnarray*}
so $T_X(n)$ is the number of legal-covers of $L_n$ with classification $X.$
\end{Definition}

\begin{figure*}[p]
$$\begin{array}{c}
\subfigure[$CC_1$ (all edges) and $T_1$ (solid edges)]{\epsfig{file=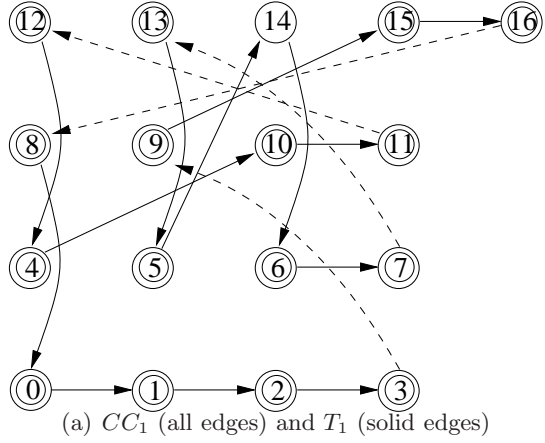, width=2.8 in}} \\[.3in]
\subfigure[$CC_2$ (all edges) and $T_2$ (solid edges)]{\epsfig{file=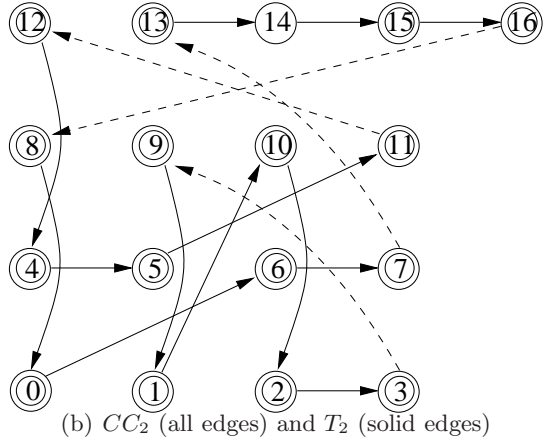, width=2.8 in}} \\[.3in]
\subfigure[$CC_3$ (all edges) and $T_3$ (solid edges)]{\epsfig{file=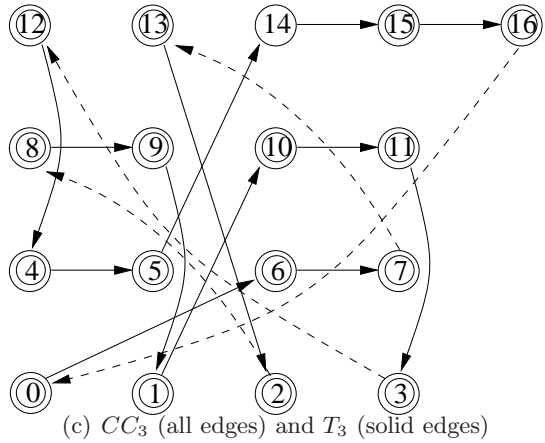, width=2.8 in}}
\end{array}$$
\caption{All of the figures are subsets of $C_{4n+1}^{1,n+2,2n+1}$. 
Solid edges are in $L_n$; dashed edges are in $\Hookn$.
The union of solid and dashed edges 
 comprise different cycle covers $CC_i$, $i=1,2,3$ in $C_{4n+1}^{1,n+2,2n+1}$.
Removing the dashed edges leaves three legal covers $T_i$, $i=1,2,3$ in $L_n$.
$C(T_1)=C(T_2)=((1,1,1,1,0,0,0,0),(0,1,0,1,0,1,0,1)).$
$C(T_3)=((0,1,1,1,0,1,0,0),(0,0,0,1,1,1,0,1)).$
}
\label{fig:non_const_cc}
\end{figure*}

\begin{figure*}[p]
$$\begin{array}{c}
\subfigure[$T_1 \cup \{(3,14),(13,20),(14,4),(19,9)\}$]{\epsfig{file=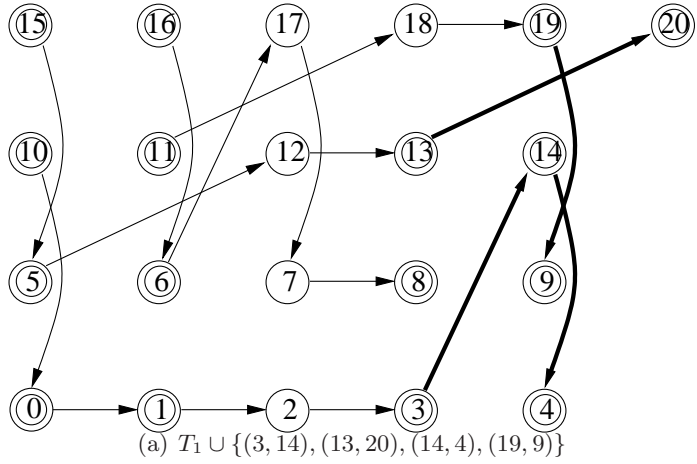, width=3.6 in}} \\[.3in]
\subfigure[$T_2 \cup \{(3,14),(13,20),(14,4),(19,9)\}$]{\epsfig{file=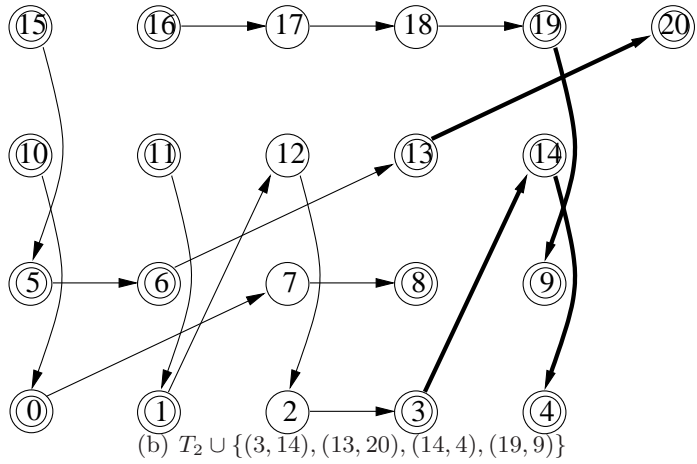, width=3.6 in}} \\[.3in]
\subfigure[$T_3 \cup \{(3,14),(13,20),(14,4),(19,9)\}$]{\epsfig{file=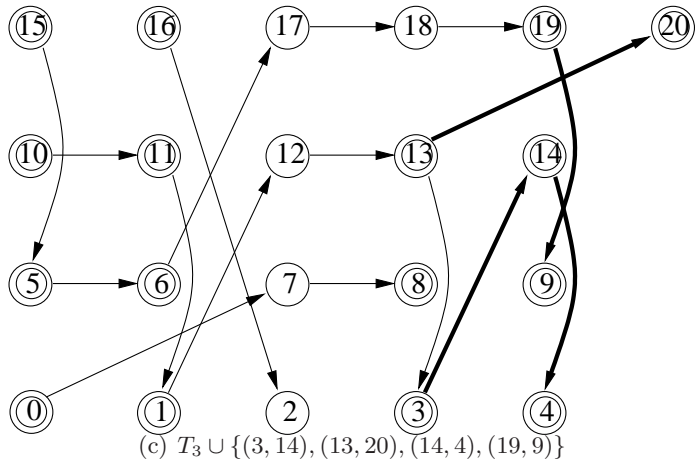, width=3.6 in}}
\end{array}$$
\caption{$n$ was increased from $4$ to $5$ and 
$S=\{(3,14),(13,20),(14,4),(19,9)\}$ was added to $T_i$, of previous figure.
Note that, $C(T_3 \cup S)=\emptyset$ since it is no longer a legal
cover (see vertex $13$).
$C(T_1 \cup S)=C(T_2 \cup S)=((1,1,1,1,0,0,0,0),(0,1,0,0,1,1,0,1))$.
}
\label{fig:non_const_class}
\end{figure*}

\begin{Lemma} See Figures \ref{fig:non_const_cc} to \ref{fig:non_const_class}. \\
 Let $ X=(L^X, R^X) \in \calp$.
Let $T_1$ be a legal cover in $L_{n_1}$ and $T_2$ be a legal cover of $L_{n_2}$, such that
$C(T_1)=C(T_2) =X$.

(a) Let   $S \subseteq \Hookn.$   Then, \vspace*{.1in}\\
\centerline{$T_1 \cup S$ is a cycle-cover of $C_{n_1}$}\\
\centerline{ {\bf iff}}\\
\centerline{ $T_2 \cup S$ is a cycle-cover of $C_{n_2}$} \\

\vspace*{.1in}

\noindent(b) 
Let  $S \subseteq \Nn.$ Then, 
$$C(T_1 \cup S) = C(T_2 \cup S).$$
That is,  either {\em both} $T_1 \cup S$ and $T_2 \cup S$ are {\em not} legal  covers
or, they are both legal covers and there is some $X'  \in \calp$ such that
$C(T_1 \cup S) = C(T_2 \cup S) = X'$
\end{Lemma}

\begin{Definition} 
For $X,X' \in \cal P,$ $S \subseteq \Hookn$ and $S' \subseteq \Nn$ set
$$
\beta_{X,S} = \left\{
\begin{array}{ll}
1 & \mbox{ if $X \cup S$ is  a cycle cover}\\
0 & \mbox{ otherwise}
\end{array}
\right.
\quad\mbox{\rm and} \quad
\alpha_{X,X',S'} = \left\{
\begin{array}{ll}
1 & \mbox{ if $C(X' \cup S')=X$}\\
0 & \mbox{ otherwise}
\end{array}
\right.
.
$$
Now set
$$
\beta_X = \sum_{S \subseteq \Hookn} \beta_{X,S}
\quad \quad 
\mbox{\rm and}
\quad \quad 
\alpha_{X,X'} = \sum_{S' \subseteq \Nn} \alpha_{X,X',S'}.
$$
\end{Definition}

Because
$NV(n)$ is no longer just the one vertex set $\{n\}$,
Lemma \ref{lem:alpha_calc} has to be replaced by
\begin{Lemma}
\label{lem:alpha_calc_non_con}
If $\alpha_{X,X',S'} =1$, then  $|S'| =|NV(n)| = p$.
\end{Lemma}
The proof is very similar to that of Lemma \ref{lem:alpha_calc}.
We now continue with 
\begin{Lemma}
$$T(n) = \sum_{X \in \calp} \beta_X T_X(n)$$
and
$$
T_X(n+1) = \sum_{X' \in \calp} \alpha_{X,X'} T_X(n).
$$
\end{Lemma}

We   reuse the concept of consistent ordering introduced in 
Definition \ref{def:consistent}.  
It is now straightforward to redo the steps of the proof of
Lemma \ref{lem:high_deg_bound} to prove
\begin{Lemma}
\label{lem:high_deg_bound_nc}
Let $A = (\alpha_{X,X'}).$ If $X \in \calp$ is ordered consistently, then there exists
an $2^{p\ems} \times 2^{p\ems}$ matrix $\bar{A}$ such that,
\[
\left(
   \begin{array}{cccc}
      \bar A & 0 & \cdots & 0 \\
      0 & \bar A & 0 & 0 \\
      \vdots & 0 & \ddots & \vdots \\
      0 & 0 & \cdots & \bar A
   \end{array}
\right)
\]
i.e. $A=\mbox{diag}(\bar A,\bar A,\ldots,\bar A)$ where $A$ contains
$2^{p\ems}$ copies of $\bar A$ on its diagonal. 
\end{Lemma}
To see that this really is a tight generalization of Lemma \ref{lem:high_deg_bound}
note that in the constant case  $p=1$ and $s=0$ and Lemma \ref{lem:high_deg_bound_nc}
then says that
the size of $\bar A$ is $2^{\bar s}$ which is exactly the result in
Lemma \ref{lem:high_deg_bound}.
 

The main difference between the constant-jump and non-constant jump case is that,
in the constant-jump case we were able,  in Lemma \ref{lem:degree_bound},  to reduce the order  of
the recurrence relation from the size of $\bar A$ to one less than the size of
$\bar A$.  This was done by using special structural properties of $\bar A$.  One of
the facts that implicitly contributed to these properties was that the {\em size} of
$NV(n)$, i.e., the number of new vertices added at each step,  was equal to one.
This is not true in the non-constant jump case and we are  therefore not able
to extend Lemma \ref{lem:degree_bound} here.
So,  the best that we can get,  from Lemma \ref{lem:high_deg_bound_nc} is that
the recurrence relation
$T(n)$ satisfies  a degree-$\left(2^{p\ems}\right)$ polynomial, an improvement of
a factor of $2^{p\ems}$ over the naive solution.

For an example of such a recurrence relation, see the second
set of graphs in Table \ref{tab:cov_results}.


\section{Variations and Extensions}
In this section we sketch some extensions to the result in the paper as well as some other uses of
the transfer matrix technique presented.  For clarity,  the results are only
shown for constant-jump circulants.  Using the techniques of  Section \ref{sec:non_constant}
 is straightforward to generalize the results
in this section to non-constant jump circulants as well.
\label{sec:Extensions}
\subsection{Weighted Circulants}
Until now we have assumed that our circulants are 0-1 matrices corresponding to being the
adjacency matrices of circulant graphs.  The permanents then counted the number of cycle covers in the corresponding circulant graphs.  
An obvious generalization is to permit the nonzero $a_i$ to be arbitrary values. 

In this case the matrix  becomes a weighted adjacency matrix.  For subsets  $T$  of the 
edges in $C_n$  let {\em weight} of $T$ be $w(T) = \prod_{(i,j) \in T} a_{i,j}.$  In this case
the permanent is the sum of the weights of all cycle covers in the corresponding $C_n$, i.e., 
$T(n) = \sum_{T \in {\cal CC}_n} w(T).$  We can modify our technique by changing the definition
of $T_x(n)$ in Definition \ref{def:class_def} to
$$T_X(n)     = \sum_{T \in {\cal L}_X(n)} w(T)$$
and the definitions of $\beta_{X,S}$ and $\alpha_{X,X',S'}$
in Definition \ref{def:alph_bet} to be
$$
\beta_{X,S} = \left\{
\begin{array}{ll}
w(S) & \mbox{ if $X \cup S$ is  a cycle cover}\\
0 & \mbox{ otherwise}
\end{array}
\right.
$$
and 
$$
\alpha_{X,X',S'} = \left\{
\begin{array}{ll}
w(S) & \mbox{ if $C(X' \cup S')=X$}\\
0 & \mbox{ otherwise}
\end{array}
\right.
.
$$

With these changes the rest of the derivations and analyses remain the same and all
of the Lemmas and proofs follow accordingly.  In particular,  we can show that
the permanent still satisfies a degree($2^{\ems} -1$) recurrence relation.

\subsection{Counting Cycles in Restricted Permutations}

In Section \ref{sec:Introduction} we discussed  how the permanent evaluates the {\em
  number} of restricted permutations using  the given jumps, i.e.,
$T(n)$ also counted  the number of permutations in  
$$
S_{n}(S) = \{ \pi \in S_{n} \,:\, \pi[i] -i \bmod (n) \in S\}.
$$
We can easily modify the transfer matrix technique to answer other questions about these permutations.
As an example,  suppose that we pick a permutation $\pi$ uniformly at random from
$S_{n}(S)$ and set $X = \mbox{\# of cycles in $\pi$}$. What can be said about
the moments of $X?$

First assume that, as previously,   $0 = s_1 < s_2 < \cdots < s_k$.
Suppose now that for cycle cover $T \in \calcc(n)$ we define
$\#_C(T)$ to be the {\em number} of cycles composing cover $T$ and set
$$TC_i(n) = \sum_{T \in \calcc(n)} \left(\#_C(T)\right)^i.$$
That is, $TC_0(n) = T(n)$ while $TC_1(n)$ is the total number of
cycles summed over all cycle-covers in $\Cn.$  Then, again by the correspondence,
we have that the moments of $X$ are given by 
$$\forall i \ge 0,\quad E(X^i) = \frac {TC_i(n)} {TC_0(n)}.$$
The interesting point is that the transfer matrix approach introduced in 
this paper can mechanically be extended 
through appropriate changes to the definition 
of $T_X(n)$ in Definition \ref{def:class_def} 
and the definitions of $\beta_{X,S}$ and $\alpha_{X,X',S'}$
in Definition \ref{def:alph_bet}, to permit showing that
for every $i,$  $TC_i(n)$ satisfies a fixed-order 
constant coefficient  recurrence relation. For given, 
$s_1,s_2,\ldots,s_k$ this permits, for example,
calculating $E(X)$ and $Var(X)$.
 
We should note that we are only saying that for the cost functions, $TC_i(n)$,
 the transfer matrix defined by Lemma \ref{lem:main_const} exists.
Lemmas \ref{lem:high_deg_bound}  and \ref{lem:degree_bound} will no longer hold, though.
So the degree of the recurrence relation will be $2^{2 \ems}$ and not $2^{\ems}-1.$.

Another complication is that in the general case   $TC_i(n)$, $i>0,$
we may  no longer assume that $0 = s_1 < s_2 < \cdots < s_k$.   Recall that we were allowed to make
this assumption when calculating the permanent ($i=0$) because   the permanent was
invariant under rotation of rows.  This is no longer true for $TC_i(n)$, $i>0.$
As an example,  consider the simple circulants $C^0_n$ (every vertex points to itself) 
and $C^1_n$(every vertex points to its neighbor).  The adjacency matrix of the first is $I_n$;
the adjacency matrix of the second $P_n^1.$  These T are rotationally
equivalent to each other.  In both cases there is only one cycle cover;  in 
$C^0_n$ it is the union of $n$ self loops;  in
$C^1_n$ the directed circle.  So,  for $C^0_n$,  $TC_1(n) = 1$ while for
$C^1_n$, $TC_1(n) = n$ and their values are different. Thus,  rotationally
equivalent circulant matrices may have different values of $TC_1(n)$.

We therefore need to modify our technique to work when $0 \not=s_1$  by
appropriately 
modifying the definition  of classifications.  The major new complication here
is that some of the edges in $\Hookn$ might be going from $L(n)$ to $R(n)$ rather 
than from $R(n)$ to $L(n)$.  
Set 
$$
\begin{array}{rclcrcl}
S^+  &=& \{s \in S \,:\, s \ge 0\},  &\quad& 	s^+ &=& \max_{s \in S^+}s,\\
S^-  &=& \{s \in S \,:\, s < 0\},			&\quad& s^-  &=& \max_{s \in S^-}|s|
		\mbox{ (if $S^- = \emptyset$ set $s^- = 0$)}
		\end{array}
$$
Now define
$$
\begin{array}{rclcrcl}
L^+(n) &=& \{0,\ldots s^+-1\}, &\quad&   R^+(n)&=& \{n-s^+,\ldots,n-1\},\\
L^-(n) &=& \{0,\ldots s^--1\},  &\quad& R^-(n)&=& \{n-s^-,\ldots,n-1\}.
\end{array}
$$
Set  $\ems = s^+ + s^-$ and 
let $\calp$ be the set of $2^{2\ems}$ tuples $(L_+,L_-, R_+,R_-)$ where 
$L_+,L_-, R_+,R_-$ are, respectively, binary $s^+,$ $s^-,$ $s^+,$ $s^-$ tuples.

\begin{Definition}
\label{def:new_legal_cov}
$T \subseteq E_L(n)$ is  a {\em legal cover of $L_n$} if 
\begin{itemize}
\item $\forall v \in V,\quad  \ID T v \le 1 \mbox{ and }\  \OD T v \le 1$.
\item $\forall v \in V - \left(L^+(n) \cup R^-(n)\right),\quad  \ID T v =1$.
\item $\forall v \in V - \left(L^-(n) \cup R^+(n)\right),\quad \OD T v =1$.
\end{itemize}
\end{Definition}

\vspace*{.1in}

Let $T$ be a {\em legal-cover} of $L_n.$ The {\em classification of $T$} will now be 
$C(T)=   (L^T_{+} ,L^T_{-}, R^T_{+},R^T_{-}) \in \calp$
where
$$
\begin{array}{lrcl}
\forall 0 \le i <  s^+,\, & L^T_+(i) &=& \ID T i\\[.04in]
                            & R^T_+(i) &=& \OD T {n-1-i},\\[.04in]
\forall 0 \le i <  s^-,\, & R^T_-(i) &=&\ID T {n-1-i},\\[.04in]
                          & L^T_-(i) &=& \OD T i.
\end{array}
$$
Not that the difference between this and the previously defined classifications was that
previously,  because $s_0 =0,$   we had $L^{-}(n) = R^{-}(n) = \emptyset.$
Given these new definitions,  we can use the same transfer matrix machinery as before
to derive recurrence relations for the $TC_i(n)$.

As an illustration recall the results from Table \ref{tab:cov_results} counting the
number of cycle covers in $C_n^{-1,0,1}$ and $C_n^{0,1,2}$.  Even though these two
graphs are {\em not isomorphic} they had the same number of cycle-covers because
the  adjacency matrix of the second  is just the adjacency matrix of the first 
with every row (cyclically) shifted over one step.  Since permanents are invariant under cyclic shifts
both matrices have the same permanent which is $\sim \phi^n$ where $\phi = (1 + \sqrt 5)/2.$

We calculated $TC_1(n)$ for both cases with the results given in  Table
\ref{table:2}.
In both cases we have that $TC_1(n) \sim c n \phi^n$. This means that if a permutation on $n$ items 
is chosen at random from the corresponding distribution then,  on average, it  will have 
$\frac {T_1(n)} {T_0(n)} \sim c n$ cycles.  
It is interesting to note that  that $c$ is  different for the two cases.

\begin{table*}[t]
$$
\begin{array}{|l|l|c| l|}
\hline
              & TC_1(n) = 3 TC_1(n-1) - TC_1(n-2) &  & \\
              & \hspace*{.6in}  -3TC_1(n-3) + TC_1(n-4) & TC_1(n) & \\
 C_n^{-1,0,1} & \hspace*{.6in}  + TC_1(n-5)     &  \sim \frac {\phi^4} {\phi^2 + \phi^4} n \phi^n & \frac {TC_1(n)} {TC_0(n)} \sim .7236n\\
   & \mbox{ initial values } 22,42,80,149,274 & \sim .7236n \phi^n & \\
   & \mbox{ for } n=4,5,6,7,8 & & \\
   \hline
              & TC_1(n) = 3 TC_1(n-1) - 6TC_1(n-3) &  & \\
              &  \hspace*{.6in} +2TC_1(n-4) + 4TC_1(n-5)  & TC_1(n) & \\
 C_n^{0,1,2} & \hspace*{.6in} - TC_1(n-6) - TC_1(n-7)    &  \sim \frac {\phi^2} {\phi^2 + \phi^4} n \phi^n
     & \frac {TC_1(n)} {TC_0(n)} \sim .2764n\\
   & \mbox{ initial values } 21, 32, 56, 93,161, 275, 475 & \sim .2764n \phi^n & \\ 
   & \mbox{ for } n=4,5,\ldots,10 & & \\
   \hline
\end{array}
$$
\caption{$TC_1(n)$ is the number of cycles summed over all cycle covers in the given graph with $n$ vertices. $TC_0(n)= T(n)$ is the number of cycle covers.  In Table \ref{tab:cov_results}
we saw that, in both cases, $TC_0(n) \sim \phi^n$  where $\phi = (1 + \sqrt 5)/2.$ }
\label{table:2}
\end{table*}

\subsection{Hamiltonian Cycles and Other Problems}
Finally,  we note that   a minor modification to the transfer-matrix technique 
permits using it to show that the number of {\em Hamiltonian Cycles}
in a directed circulant graph $C_n$ also satisfies a constant-coefficient
recurrence relation in $n$. This fact was previously known for {\em undirected} 
circulant graphs \cite{GoLe04,GoLeWa04} but doesn't seem to have been
known for directed circulants, with the exception of the special
case of in(out)-degree 2 circulants \cite{YaBuCeWo97}, also known as 
{\em two-stripe} circulants.

Again,  as when calculating $TC_1(n)$  in the previous subsection,  we may no longer assume
that $s_1=0.$  We reuse the definitions of
$L^+(n),L^{-}(n), R^+(n),R^-(n)$ introduced above and define a 
\begin{Definition}
\label{def:new_legal_cov_HC}
$T \subseteq E_L(n)$ is  a {\em legal tour  of $L_n$} if 
(i) $T$ is a Hamiltonian Cycle of $L_n$ or (ii)
\begin{itemize}
\item $\forall v \in V,\quad  \ID T v \le 1 \mbox{ and }\  \OD T v \le 1$.
\item $\forall v \in V - \left(L^+(n) \cup R^-(n)\right),\quad  \ID T v =1$.
\item $\forall v \in V - \left(L^-(n) \cup R^+(n)\right),\quad \OD T v =1$.
\item $T$ contains no cycles
\end{itemize}
\end{Definition}
Note that if $T$ is legal and is not Hamiltonian, then $T$ is composed of
paths in which (i) the start of each path is in $R^+(n) \cup L^-(n)$ (ii) 
the end of each path is in $L^+(n) \cup R^-(n)$ and (iii) every vertex is on 
exactly one path (if  a vertex $v$ is isolated we consider it to be lying
on a zero-length path that starts and ends at $v$).  The {\em classification} of
$T$ will then be the union of the (start,end) pairs describing the starting and ending
points of each path.  The number of such classifications is finite.  Furthermore,  the 
classification of $T \cup S$ where $T$ is a legal tour and 
$S \subseteq \Nn$ or $S\subseteq \Hookn$ depends only upon the classification of $T$ 
and the edges in $S$.  We can therefore use the method described in this paper to show
that the number of Hamiltonian cycles in $C_n$ satisfies a recurrence relation.

We conclude by noting that there is nothing particularly special about 
Hamiltonian Cycles and that the technique will enable counting many other
structures in directed circulant graphs  as well.  As an example,  it is not too difficult to
modify the method to show that  the number of Eulerian Tours in such graphs also satisfies a constant-coefficient
recurrence relation in $n$.

%
%

\section{Conclusion}
In this paper we showed a new derivation of Minc's result  \cite{Minc78,Minc87b}) that the
permanent of parametrized circulant matrices satisfies a recurrence relation.  Instead 
of being algebraic our new technique
was combinatorial.  We took advantage of the fact that permanents of 0/1 matrices
count the number of directed cycle covers in the matrix associated
with the graph to transform the problem into a counting one.  We were then able to
decompose circulant matrices in such a way as to allow the use of the transfer matrix
method to count the number of  cycle covers.  Finally,  we were able to show
that the transfer matrix was block diagonal,  with all blocks being copies of each other,
reducing the order of the characteristic polynomial of the transfer matrix (and thus of
the corresponding recurrence relation for the permanents). 

A benefit of this new derivation is that it easily extends to the analysis of
non-constant (linear) jump circulants,  something that the 
original Minc result could not handle.  It also permits counting many other properties
of circulant graphs,  e.g., the number of Hamiltonian cycles.

{
\bibliographystyle{plain}
\bibliography{bibtex}
}


\appendix



\section{A Worked Example for  $C_n^{0,1,2}$}
\label{App:worked example}

In Sections \ref{sec:decomposition} and \ref{sec:Minc} we derived that 
 $T(n)$, the number of cycle covers in
$C_n^{0,1,2}$, satisfies
$$\forall n \ge 2\ems,\quad 
T(n) = \beta\,   \bar T (n)
\quad \mbox{and} \quad 
\bar T(n+1) =A \, \bar T(n)
$$
where 
$\beta = (\beta_X)_{X \in \calp}$ and $A = (\alpha_{X,X'})_{X,X' \in \calp}$.

\vspace*{.1in}

For $C_n^{0,1,2}$,  $\ems =  2.$
Definition \ref{def:class_def} then says that
 every $X \in \calp$ is in the form $X=(L^X,R^X)$
where $L^X,R^X \in \{0,1\}^2$.  We can therefore
represent every $X$ by a four-bit binary vector in which the first two bits represent
$L^X$ and the last  two $R^X$; there are $16$ such $X \in \calp$.

Ordering the $X$  lexicographically we calculate that $\beta$ is
\[\small
\left(
   \begin{array}{cccccccccccccccc}
       1 & 0 & 0 & 0 & 0 & 1 & 0 & 0 & 0 & 1 & 1 & 0 & 0 & 0 & 0 & 1
   \end{array}
\right),
\]
$\bar T(4)$ is 
\[\small  \left(
   \begin{array}{cccccccccccccccc}
       1 & 0 & 0 & 0 & 0 & 2 & 1 & 0 & 0 & 3 & 2 & 0 & 0 & 0 & 0 & 1
   \end{array}
\right)^t
\]
(where the $t$ denotes taking the transpose), and
{{\em Transfer matrix} $A$ is 
\[
\left(
   \begin{array}{rrrr|rrrr|rrrr|rrrr}
      1 & 0 & 0 & 0 & 0 & 0 & 0 & 0 & 0 & 0 & 0 & 0 & 0 & 0 & 0 & 0 \\
      0 & 1 & 1 & 0 & 0 & 0 & 0 & 0 & 0 & 0 & 0 & 0 & 0 & 0 & 0 & 0 \\
      0 & 1 & 0 & 0 & 0 & 0 & 0 & 0 & 0 & 0 & 0 & 0 & 0 & 0 & 0 & 0 \\
      0 & 0 & 0 & 1 & 0 & 0 & 0 & 0 & 0 & 0 & 0 & 0 & 0 & 0 & 0 & 0 \\ \hline
      0 & 0 & 0 & 0 & 1 & 0 & 0 & 0 & 0 & 0 & 0 & 0 & 0 & 0 & 0 & 0 \\
      0 & 0 & 0 & 0 & 0 & 1 & 1 & 0 & 0 & 0 & 0 & 0 & 0 & 0 & 0 & 0 \\
      0 & 0 & 0 & 0 & 0 & 1 & 0 & 0 & 0 & 0 & 0 & 0 & 0 & 0 & 0 & 0 \\
      0 & 0 & 0 & 0 & 0 & 0 & 0 & 1 & 0 & 0 & 0 & 0 & 0 & 0 & 0 & 0 \\ \hline
      0 & 0 & 0 & 0 & 0 & 0 & 0 & 0 & 1 & 0 & 0 & 0 & 0 & 0 & 0 & 0 \\
      0 & 0 & 0 & 0 & 0 & 0 & 0 & 0 & 0 & 1 & 1 & 0 & 0 & 0 & 0 & 0 \\
      0 & 0 & 0 & 0 & 0 & 0 & 0 & 0 & 0 & 1 & 0 & 0 & 0 & 0 & 0 & 0 \\
      0 & 0 & 0 & 0 & 0 & 0 & 0 & 0 & 0 & 0 & 0 & 1 & 0 & 0 & 0 & 0 \\ \hline
      0 & 0 & 0 & 0 & 0 & 0 & 0 & 0 & 0 & 0 & 0 & 0 & 1 & 0 & 0 & 0 \\
      0 & 0 & 0 & 0 & 0 & 0 & 0 & 0 & 0 & 0 & 0 & 0 & 0 & 1 & 1 & 0 \\
      0 & 0 & 0 & 0 & 0 & 0 & 0 & 0 & 0 & 0 & 0 & 0 & 0 & 1 & 0 & 0 \\
      0 & 0 & 0 & 0 & 0 & 0 & 0 & 0 & 0 & 0 & 0 & 0 & 0 & 0 & 0 & 1
   \end{array}
\right)
\]
}
The Lexicographic ordering is consistent so,
as predicted by Lemma \ref{lem:degree_bound},  $A$
is  partitioned into $16$ $ 4 \times 4$ blocks where
all but the diagonal blocks are $0$ and all of the diagonal blocks are equal to some
$4 \times 4$ matrix $\bar A$ which in this case is 
$$\bar A = 
\left(
   \begin{array}{rrrr}
     1 & 0 & 0 & 0\\
     0 & 1 & 1 & 0 \\
     0 & 1 & 0 & 0\\
     0 & 0 & 0 & 1
   \end{array}
\right).
\]
Note that the lexicographic ordering on four-bit vectors also has the property that
$$(0,0) < (0,1) < (1,0) < (1,1).$$
This means that if $X_1=(L_1,X_1)$, $X_2=(L_2,R_2)$ and $L_1=L_2$ then if the number of 
'$0$'s in $R_1$ is less than the number of '$0$'s in $R_2$ then $X_1 < X_2$.
This satisfies the conditions of the ordering used in the proof of 
Lemma \ref{lem:degree_bound} which then implies that $\bar A$ should be
in the form
$$\bar A = 
\left(
   \begin{array}{rrr}
     B_0 & 0   & 0\\
     0   & B_1 & 1 \\
     0   & 1   & B_2\\
   \end{array}
\right)
\]
where $B_i$ is a ${2 \choose i} \times {2 \choose i}$ matrix.  We do observe this
behavior with  $B_0 = B_2 = (1)$ and
$B_1 = 
\left(
   \begin{array}{rr}
     1 & 1\\
     1 & 0
   \end{array}
\right).
$
The characteristic polynomial of $B_0$ and $B_2$ is $P_0(x)=x-1$.
The characteristic polynomial of $B_1$ is $P_1(x) = x^2-x-1$.

This implies that 
$$Q(x) = P_1(x) P_0(x)  =(x^2-x-1)(x-1) = x^3 - 2x^2 +1$$
annihilates $A$.

Working through the details we can then solve to find that,
for $C_n^{0,1,2},$ 
$T(n) =  2T(n-1) - T(n-3)$ with initial values 
$T(4) = 9,$ $T(5) = 13,$ and $T(6) = 12.$

\section{Proofs of Lemmas \ref{lem:L_n} and \ref{lem:RecC_non_const} }
\label{sec:appB}

{\par\noindent\underline{\bf Proof  of Lemma \ref{lem:L_n}:}}\\[.1in]
We only prove the $\Rightarrow$ part. The reverse direction can be proved by the
same argument.

(a) If $((u_1,v_1),(u_2,v_2)) \in {\widehat E}_L(n)$, there exist $i$  such that 
\begin{eqnarray}
\nonumber f(n;u_2,v_2)-f(n;u_1,v_1) &=& p_i n + s_i \bmod{(pn+s)} \\
\nonumber &\mbox{and}&\\
u_2-u_1 &=& p_i \bmod{p} .\nonumber
\end{eqnarray}

If $f(n;u_2,v_2)\ge f(n;u_1,v_1)$, then 
$$ f(n;u_2,v_2)-f(n;u_1,v_1) = p_i n + s_i \mbox{ and } u_2-u_1 = p_i. $$
When $n$ is increased to $n+1$,
$$
\begin{array}{cl}
  & f(n+1;u_2,v_2)-f(n+1;u_1,v_1) \\
= & u_2(n+1)+v_2-u_1(n+1)-v_1 \\
= & (u_2n+v_2-u_1n-v_1)+(u_2-u_1) \\
= & p_i (n+1) + s_i
\end{array}
$$

(b) If $f(n;u_2,v_2)<f(n;u_1,v_1)$, then
$$ pn+s+f(n;u_2,v_2)-f(n;u_1,v_1) = p_i n + s_i \mbox{ and } p+u_2-u_1 = p_i. $$
When $n$ is increased to $n+1$,
$$
\begin{array}{cl}
  & p(n+1)+s+f(n+1;u_2,v_2)-f(n+1;u_1,v_1) \\
= & p(n+1)+s+u_2(n+1)+v_2-u_1(n+1)-v_1 \\
= & (pn+s+u_2n+v_2-u_1n-v_1)+(p+u_2-u_1) \\
= & p_i (n+1) + s_i
\end{array}
$$
Therefore, in both cases, $((u_1,v_1),(u_2,v_2)) \in {\widehat E}_L(n+1)$.
{{\hfill{\Large$\Box$}}\smallskip}

\bigskip

{\par\noindent\underline{\bf Proof of Lemma \ref{lem:RecC_non_const}:}}\\[.1in]
We split the proof into two parts.

\noindent  (a) \underline{$\Hookn \subseteq R(n) \times L(n):$}

Let $e=((u_1,v_1),(u_2,v_2))$ be an edge in ${\widehat E}_C(n)$ associated with the jump $p_in+s_i$.
Note that $e \in \Hookn$ if and only if
$$f(n;u_2,v_2)-f(n;u_1,v_1)=p_in+s_i \bmod{pn+s}\mbox{ and }u_2-u_1 \not= p_i \bmod{p}$$
There are two cases: \\

\noindent (i) $f(n;u_1,v_1) \le f(n;u_2,v_2)$.
\begin{eqnarray*}
&\Rightarrow & u_2n+v_2 = u_1n+v_1 + p_in+s_i \\
&\Rightarrow &(u_2-u_1-p_i)n = s_i+v_1-v_2
\end{eqnarray*}

$u_2-u_1 \neq p_i \bmod p$ implies $u_2-u_1-p_i\neq 0.$ On the other hand, 
$$ -2n < s_i+v_1-v_2 < 2n,$$ hence $u_2-u_1-p_i = \pm 1.$

If $u_2-u_1-p_i = 1,$ 
then $s_i+v_1-v_2 = n,$ so $v_1 \geq n-s_i \geq n - \ems$.
Since,  by definition,
$v_1 \leq n-1$, we also have
 $v_2 \leq s_i-1.$
Furthermore, $u_1 = u_2 -p_i -1 < p-1$.  
 Hence
$e \in R(n) \times L(n).$

If $u_2-u_1-p_i=-1,$  this implies $s_i+v_1-v_2+n = 0.$ We 
then have $v_2 \geq n+s_i \geq n+s,$
which is not possible since it's outside the range of $v_2$.

\vspace{0.1in}
\noindent (ii) $f(n;u_1,v_1) > f(n;u_2,v_2)$.
\begin{eqnarray*}
&\Rightarrow& u_2n+v_2+pn+s = u_1n+v_1 + p_in+s_i\\
& \Rightarrow &(u_2-u_1+p-p_i)n = s_i-s -v_2 +v_1.
 \end{eqnarray*}

Similar to the previous case:
$$ -2n < s_i -s -v_2 +v_1 < 2n,$$
and  $u_2-u_1 \neq p_i \bmod p$ implies $ u_2-u_1+p-p_i = \pm 1.$

If $u_2-u_1+p-p_i = 1,$ then  $s_i-s-v_2+v_1 =n.$ Thus $v_1 \geq n+s-s_i,$ and $v_2
\leq s_i -1.$ Hence $e \in R(n) \times L(n).$

If $u_2 -u_1+p-p_i = -1,$ this implies  $s_i-s-v_2+v_1 +n =0.$  We then have
 $u_2 = u_1 +p_i-p-1
\leq p_i -2 < p-1$ which implies  $v_2 \leq n-1.$ However this results in   $v_2 = s_i -s +n +v_1
\geq n+s_i-s \geq n$ which is  not possible since it's outside the range of $v_2$.


\vspace{0.2in}
Therefore, $\Hookn \subseteq \left(R(n) \times L(n)\right).$

\vspace{0.3in}

\noindent (b) \underline{$\Nn \subseteq \left(R(n) \times NV(n) \right) \cup \left( NV(n)
  \times NV(n) \right):$}  

From Lemma~\ref{lem:L_n}
$$\Nn \subseteq ({\widehat V}(n) \times NV(n)) \cup (NV(n) \times {\widehat
  V}(n)) \cup (NV(n) \times NV(n)).$$  
Let  $e=((u_1,v_1),(u_2,v_2)) \in \Nn$ be 
associated with jump $p_i(n+1)+s_i$. 

First assume $(u_1,v_1)\in {\widehat V}(n)$ 
Consider the edge $e'$  starting with $(u_1,v_1)$  associated with 
jump $p_in+s_i$ in ${\widehat E}_C(n)$, i.e.,  in circulant graph $\Cn$ and not  lattice graph 
$\Lnpo$. Then
$$e \in \Nn \Rightarrow  e \not\in {\widehat E}_L(n)
\Rightarrow e' \in \Hookn.$$
So, from part (a), $(u_1,v_1)$ is in $R(n).$ Because $(u_2,v_2) \in NV(n),$
$((u_1,v_1),(u_2,v_2)) \in R(n) \times NV(n)$. 

Now assume  that $(u_2,v_2)\in {\widehat V}(n)$.
Consider the edge $e'$ ending with $(u_2,v_2)$  associated with 
jump $p_in+s_i$ in ${\widehat E}_C(n)$. Again
$$e \in \Nn \Rightarrow  e \not\in {\widehat E}_L(n) \Rightarrow e \in \Hookn.$$
So, from part (a), $(u_1,v_1) \in L(n).$ However $L(n) \cap NV(n) = \emptyset,$ so such  a $(u_1,v_1)$
does not exist.
{{\hfill{\Large$\Box$}}\smallskip}

\end{document}